\documentclass[11pt]{amsart}
\usepackage{amscd,amssymb,verbatim,graphicx}
\usepackage[small,nohug,heads=vee]{diagrams}
\diagramstyle[labelstyle=\scriptstyle]
\let\oldmarginpar\marginpar
\renewcommand\marginpar[1]{\-\oldmarginpar[\raggedleft\footnotesize #1]%
	{\raggedright\footnotesize #1}}

\newskip\stdskip                      % standard vertical space
\newcommand{\Def}{\operatorname{Def}}

\newcommand{\Art}{\operatorname{Art}}
\newcommand{\del}{\partial}

\renewcommand{\mod}{\operatorname{mod}}

\newcommand{\und}{\underline}

\newcommand{\OO}{\mathcal{O}}

\newcommand{\bL}{{\bf L}}

\newcommand{\G}{\mathbb{G}}

\newcommand{\hra}{\hookrightarrow}

\newcommand{\CC}{\mathcal{C}}
\newcommand{\UU}{\mathcal{U}}

\newcommand{\Spec}{\operatorname{Spec}}

\newcommand{\Proj}{\operatorname{Proj}}
\renewcommand{\P}{{\Bbb P}}

\newcommand{\si}{\sigma}

\newcommand{\Pic}{\operatorname{Pic}}

\newcommand{\de}{\delta}

\newcommand{\im}{\operatorname{im}}

\numberwithin{equation}{subsection}
\newcommand{\Per}{\operatorname{Per}}

\newcommand{\A}{{\Bbb A}}

\newtheorem{thm}{Theorem}[subsection]
\newtheorem{prop}[thm]{Proposition}
\newtheorem{lem}[thm]{Lemma}
\newtheorem{cor}[thm]{Corollary}
{  \theoremstyle{definition}
\newtheorem{defi}[thm]{Definition}
\newtheorem{ex}[thm]{Example}

\newtheorem{rem}[thm]{Remark}

}

\newcommand{\Pf}{\noindent {\it Proof}}

\newcommand{\ov}{\overline}
\newcommand{\we}{\wedge}

\newcommand{\MM}{\mathcal{M}}
\newcommand{\TT}{\mathcal{T}}

\newcommand{\Hom}{\operatorname{Hom}}

\newcommand{\Ext}{\operatorname{Ext}}

\newcommand{\Res}{\operatorname{Res}}

\newcommand{\om}{\omega}

\newcommand{\la}{\lambda}

\newcommand{\C}{{\Bbb C}}

\newcommand{\Z}{{\Bbb Z}}
\newcommand{\Q}{{\Bbb Q}}

\newcommand{\wt}{\widetilde}
\newcommand{\ot}{\otimes}

\newcommand{\sub}{\subset}
\newcommand{\ed}{\qed\vspace{3mm}}

\newcommand{\cha}{{\bf char}}

\title{A modular compactification of $\MM_{1,n}$ from $A_\infty$-structures}
\author[Yank{\i} Lekili and Alexander Polishchuk]{Yank{\i} Lekili \\ Alexander Polishchuk}
\thanks{Y.L. is partially supported by a Royal Society Fellowship and a Marie Curie grant. A.P. is supported in part by the NSF grant DMS-1400390}

\address{King's College London}
\address{University of Oregon}

\begin{document}

\begin{abstract}
We show that a certain moduli space of minimal $A_\infty$-structures coincides with the modular
compactification $\ov{\MM}_{1,n}(n-1)$ of $\MM_{1,n}$ constructed by Smyth in \cite{Smyth-I}.
In addition, we describe these moduli spaces and the universal curves over them by explicit equations,
prove that they are normal and Gorenstein, show that their Picard groups have no
torsion and that they have rational singularities if and only if $n\le 11$.
\end{abstract}

\maketitle

\centerline{\sc Introduction}

\quad\quad One of the motivations of the present work is to show that a study of the
derived categories of coherent sheaves on such basic varieties as algebraic
curves can uncover interesting geometry, including some aspects of the moduli
spaces.  The idea to study algebraic varieties via their derived categories,
which has been around for a while (see \cite{BO}), recently got more focus and
motivation coming from the homological mirror symmetry. In particular, it
became clear that it is important to take into account the dg-enhancement, or
the corresponding $A_\infty$-structure obtained by homological perturbation (in
other words, one has to keep track of the higher Massey operations). Namely, if
one takes a generator $G$ of the derived category, then the corresponding
$\Ext$-algebra $\Ext^*(G,G)$ has a structure of an $A_\infty$-algebra, from
which the derived category can be recovered.  This raises a natural question:
what kind of $A_\infty$-algebras are obtained in this way, possibly for some
specially chosen generators $G$. For example, for a smooth projective curve $C$
we can take as a generator of the derived category the object
\begin{equation}\label{generator-C-points-eq}
	G(C,p_1,\ldots,p_n)=\OO_C\oplus\OO_{p_1}\oplus\ldots\oplus\OO_{p_n},
\end{equation} where $p_1,\ldots,p_n$ are distinct points of $C$. In the case
$n=g$, the genus of $C$, the resulting $A_\infty$-algebras were studied in
\cite{LP}, \cite{LP2} (for $g=1$) and \cite{P-ainf} (in general). The case of genus
$0$ curves was also studied in \cite{P-ainf}. In this paper we consider
the case when $C$ is of genus $1$ and $n$ is arbitrary.

Note that to recover the derived category $D^b(C)$ from the $A_\infty$-algebra $E_{C,p_1,\ldots,p_n}$
associated with the generator \eqref{generator-C-points-eq}
one only needs the category of $A_\infty$-modules over $E_{C,p_1,\ldots,p_n}$.
On the other hand, this $A_\infty$-algebra itself, viewed up to a gauge equivalence,
carries more information: in fact, one can recover the pointed curve $(C,p_1,\ldots,p_n)$ from it.
Moreover, in some situations one gets 
an equivalence between the appropriate moduli spaces of curves 
and moduli spaces of $A_\infty$-algebras. 

In order to get such an equivalence, one has to allow curves to be singular but 
also impose some restrictions on $(C,p_1,\ldots,p_n)$ guaranteeing that
the associative algebra structure on $\Ext^*(G,G)$ for $G=G(C,p_1,\ldots,p_n)$
is independent of the pointed curve (whereas higher products do depend on it), and that $G$ is indeed a generator
of the perfect derived category of $C$. The latter property is equivalent to the ampleness of $\OO_C(p_1+\ldots+p_n)$.

The above program was implemented in \cite{LP} and \cite{P-ainf} for the case $n=g$.
In this paper we study a similar equivalence between moduli of curves and $A_\infty$-structures in 
the case of curves of arithmetic genus one with $n>1$ (smooth) marked points.
In this case the algebra $\Ext^*(G,G)$ does not depend on a curve provided one has $H^1(C,\OO(p_i))=0$ for
each of the marked points. In addition, we require $\OO_C(p_1+\ldots+p_n)$ to be ample.
We call the resulting moduli stack $\UU_{1,n}^{sns}$ (``sns" stands for ``strongly non-special", since
each $p_i$ defines a non-special divisor). 

Note that the relevance of our work to symplectic geometry is due to the fact that 
Fukaya category of $n$-marked (symplectic) torus gives rise to an $A_\infty$-structure in
$\UU_{1,n}^{sns}$. An attempt to directly compute this $A_\infty$-algebra runs
into the well-known transversality problems: the constant maps contribute non-trivially to
higher products (cf. \cite{LP}, \cite{LP2} for $n=1$). To deal with this problem, one has to device
a consistent set of perturbations which makes the computations hard (as one has to solve infinitely many PDEs in a consistent way). Our result tells that once
we know that the cohomology algebra is isomorphic to $\Ext^*(G,G)$ (which is easy to check), then we
know that at the chain level the $A_\infty$-algebra arising from the Fukaya category corresponds to
one of the curves in the moduli space $\UU_{1,n}^{sns}$. From this characterization it follows that
whatever one wants to compute for a given $A_\infty$-structure in our moduli space, in particular
the one coming from the Fukaya category, we can do so using the {\it commutative} model given by the derived category of the corresponding curve. In the follow-up work \cite{LPol}, extending the ideas of
\cite{LP}, we will use this to
establish a very precise form of homological mirror symmetry for $n$-punctured tori which is valid over $\Z$. 

The paper consists of two parts. The first is a purely algebro-geometric study of the moduli stacks $\UU_{1,n}^{sns}$
(without any reference to $A_\infty$-structures). Here our main result identifies $\UU_{1,n}^{sns}\setminus \{ C_{1,n}\}$, 
where $C_{1,n}$ is the {\it elliptic $n$-fold curve} (a certain generalization of the cuspidal cubic
curve, see \ref{Smyth-sec}), with one of the alternative compactifications of $\MM_{1,n}$ constructed and studied by Smyth in \cite{Smyth-I, Smyth-II}.
Recall that for each $m$, $1\le m<n$, Smyth constructs the moduli stack $\ov{\MM}_{1,n}(m)$ of {\it $m$-stable} 
$n$-pointed curves. The definition of $m$-stability involves restricting the types of singularities a curve can have (see Def. \ref{mstab}). 

Smyth proves that these are proper irreducible Deligne-Mumford stacks over $\Spec(\Z[1/6])$, that the
corresponding coarse moduli spaces are projective and for $m\le 10$ coincide with certain log
canonical models of the Deligne-Mumford compactification $\ov{M}_{1,n}$
proposed by Hassett and Keel (\cite{Smyth-II}, \cite{hassett}). Our interest in
the current study is the case $m=n-1$, and we consider a slight modification of the Smyth's moduli stack, 
denoted by $\ov{\MM}_{1,n}^\infty$ (see \ref{modified} for precise definition),
which is a proper algebraic stack over $\Spec(\Z)$ (resp., $\Spec(\Z[1/2])$ for $n=2$; resp., $\Spec(\Z[1/6])$ for $n=1$),
such that
$\ov{\MM}_{1,n}(n-1)=\ov{\MM}_{1,n}^\infty\times\Spec(\Z[1/6])$.

\medskip

\noindent
{\bf Theorem A}[=Thm.\ \ref{moduli-curves-thm}+Thm.\ \ref{curve-moduli-thm}]. {\it Let us work over $\Spec(\Z)$ for $n\ge 3$, over $\Spec(\Z[1/2])$ for $n=2$, and over
$\Spec(\Z[1/6])$ for $n=1$. Let $\wt{\UU}^{sns}_{1,n}\to \UU^{sns}_{1,n}$ be the $\G_m$-torsor associated
with the standard line bundle $\la$. Then $\wt{\UU}^{sns}_{1,n}$ is an affine scheme of finite type, defined
by explicit equations. We have a natural equivalence
$$(\wt{\UU}_{1,n}^{sns}\setminus \{ C_{1,n} \} )/\G_m\simeq \ov{\MM}_{1,n}^\infty.$$
}

As a consequence of this approach to $\ov{\MM}_{1,n}^\infty$ we are able to
establish some additional results about geometry of this moduli space. For
example, we show that for $n\ge 5$ it is a projective scheme given by explicit quadratic equations. We prove that
$\ov{\MM}_{1,n}^\infty$ is smooth for $n \leq 6$ and has rational singularities
if and only if $n \leq 11$. For small values of $n$, we identify the
moduli space explicitly. For example, we have an amusing identification 
\[ \ov{\MM}_{1,6}^\infty \cong Gr(2,5), \] 
where $Gr(2,5)$ stands for the
Grassmannian of $2$-planes in the $5$-space. This extends the well known presentation
of elliptic normal curves of degree $5$ as linear sections of $Gr(2,5)$ (see \cite{Fisher}) to
singular curves (see Corollary \ref{linear-section-cor}).

A key structural result that we repeatedly use is the identification of the
universal curve over $\ov{\MM}_{1,n-1}^\infty$ with a blow-up of
$\ov{\MM}_{1,n}^\infty$ at $n-1$ special points, in such a way that the exceptional divisors of the blow-up
correspond to the universal marked points (see Proposition \ref{blow-up-prop}). This also allows us to deduce that
$\ov{\MM}_{1,n}^\infty$ is normal and Gorenstein. In particular, we deduce that Smyth's moduli stacks
$\ov{\MM}_{1,n}(m)$ are normal and Gorenstein, which simplifies some statements
in \cite{Smyth-II} formulated using the normalizations. 
In addition, we prove that the Picard group of $\ov{\MM}_{1,n}^\infty$ is freely generated by $\la$
(this was known rationally; we show that the Picard group has no torsion).

In the second part of the paper we identify the moduli stacks $\UU_{1,n}^{sns}$ with an
appropriate moduli of $A_\infty$-structures. Namely, we consider minimal $A_\infty$-algebra structures on the graded algebra $E_{1,n}=\Ext^*(G,G)$ for $G=G(C,p_1,\ldots,p_n)$ given by \eqref{generator-C-points-eq}.
As we mentioned above, up to an isomorphism this algebra does not depend on $(C,p_1,\ldots,p_n)$. Passing to the
$\G_m$-torsor $\wt{\UU}_{1,n}^{sns}$ allows to fix such an isomorphism canonically.
Working over a field $k$, 
we prove using \cite[Cor.\ 4.2.5]{P-ainf} that the functor of minimal $A_\infty$-structures on $E_{1,n}$, viewed up to a gauge
equivalence, is representable by an affine scheme $\MM_\infty(E_{1,n})$ of finite type over $k$.

\medskip

\noindent
{\bf Theorem B}[=Thm.\ \ref{ainf-thm}]. {\it Assume that $\cha(k)\neq 2$ if $n=2$ (resp., $\cha(k)\neq 2,3$ if $n=1$).
The affine scheme $\wt{\UU}^{sns}_{1,n}\times \Spec(k)$ is naturally isomorphic to the moduli scheme 
$\MM_\infty(E_{1,n})$
of minimal $A_\infty$-structures on the algebra $E_{1,n}$ up to a gauge equivalence.}

\medskip

The proof of Theorem B
follows the same strategy as in \cite{P-ainf}: we reduce this to the study of the deformation functors around the most
singular point of $\UU_{1,n}^{sns}$ which corresponds to the trivial $A_\infty$-structure. Thus, using the tools developed
in \cite{P-ainf} we reduce the proof to studying the deformation functor of the elliptic $n$-fold curve.
Similar approach is used by one of us in the follow-up paper \cite{P-more-pts} on the moduli of curves of genus $g$
with $n\ge g$ marked points forming a non-special divisor. However, the latter moduli stack in the case $g=1$ is larger than 
$\UU_{1,n}^{sns}$ (which it contains  as an open substack), so its identification with moduli of $A_\infty$-structures requires the stronger assumption that the characteristic of the ground field is $\neq 2$. 

On the other hand, the purely algebro-geometric side of our study 
is continued in \cite{P-krich}, \cite{P-more-pts} and \cite{P-g1}. Still, our understanding of the geometric properties
of the moduli spaces $\wt{\UU}_{1,n}^{sns}$ is much better than that of other similar moduli spaces. 

\begin{section}{Curves of arithmetic genus $1$ with $n$ marked points}
	
\begin{subsection}{Normal forms of pointed curves of arithmetic genus $1$}\label{normal-forms-sec}
\

\quad\quad Let $k$ be an algebraically closed field and let $C$ be a reduced, connected projective 
curve over $k$ of arithmetic
genus 1 with $n$ distinct smooth marked points: $p_1,\ldots,p_n$. When $C$ is
smooth, $(C,p_1,\ldots,p_n)$ defines a point of the moduli stack
$\mathcal{M}_{1,n}$. To compactify $\mathcal{M}_{1,n}$ one has to allow $C$ to be singular. More precisely,
we are interested in {\it modular} compactifications of $\mathcal{M}_{1,n}$ obtained
by specifying a deformation open class of curves satisfying the unique limit property in families (see \cite{smyth}).
For example, the Deligne-Mumford moduli space of
stable curves is a well known modular compactification for which $C$ is
allowed to have nodal singularities and $\omega_C (p_1 + \ldots + p_n)$ is required to be ample (\cite{DM}).  
Smyth's moduli spaces $\ov{\MM}_{1,n}(m)$ parametrize certain curves that are allowed to have elliptic $l$-fold points
with $l\le m$, in addition to nodes (see Section \ref{Smyth-sec} below).

In our approach we start by considering
the following requirements on $(C, p_1, \ldots, p_n)$: 
\begin{enumerate} 
	\item \label{cond1} $h^0(\mathcal{O}_C (p_i)) = 1$ for all $i$.  
	\item \label{cond2} $\mathcal{O}_C (p_1+\ldots+p_n)$ is ample.  
\end{enumerate} 		
We will see that the type of singularities of $C$ will be determined by these conditions \emph{a posteriori}.

We will first follow a pedestrian approach in constructing the moduli space of
arithmetic genus 1 curves $(C,p_1,\ldots,p_n)$ which satisfy the conditions
(\ref{cond1}) and (\ref{cond2}). Recall that a reduced curve is automatically
Cohen-Macaulay, so we have a relative dualizing sheaf $\omega_C$ on $C$. 
The condition \eqref{cond1} has a simple interpretation that we will use repeatedly.

\begin{lem}\label{diff-restriction-lem} 
(i) Let $C$ be a reduced connected projective curve of arithmetic genus $1$, $p$ a smooth point on $C$.
Then $h^0(\OO_C(p))=1$ if and only if $h^1(\OO_C(p))=0$ if and only if the restriction map
$$H^0(C,\omega_C)\to \omega_C|_p$$
is an isomorphism.

\noindent (ii) Let $\pi:C\to S$ be a flat projective morphism of relative dimension $1$, with reduced connected geometric fibers
of arithmetic genus $1$, and let
$p:S\to C$ be a section such that $\pi$ is smooth near $p$. 
Then the following conditions are equivalent:

\noindent (1) the natural map $\OO_S\to\pi_*(\OO_C(p(S)))$ is an isomorphism;

\noindent (2) $R^1\pi_*(\OO_C(p(S)))=0$;

\noindent (3) the natural map $\pi_*(\om_{C/S})\to p^*\om_{C/S}$ is an isomorphism.
\end{lem}

\Pf . (i) One uses the fact that $\chi(\OO_C(p))=1$ (by Riemann-Roch theorem)
and $h^1(\OO_C(p))=h^0(\om_C(-p))$ (by Grothendieck-Serre duality).  Since
$H^0(C,\om_C)$ is $1$-dimensional, the restriction map to $\omega_C|_p$ is an
isomorphism if and only if $h^0(\om_C(-p(S)))=0$.

\noindent (ii) The exact sequence
$$0\to \om_C(-p(S))\to \om_{C/S}\to p_*p^*\om_{C/S}\to 0$$
shows that (3) implies that the map 
$$R^1\pi_*(\om_C(-p(S)))\to R^1\pi_*(\om_{C/S})$$ 
is an isomorphism. By Grothendieck duality, this is equivalent to (1). Conversely, if (1) holds then
we get that the morphism of line bundles
$\pi_*\om_{C/S}\to p^*\om_{C/S}$ is surjective, hence, an isomorphism. Thus, (1) and (3) are equivalent.
The same exact sequence shows that (3) implies the vanishing of $\pi_*(\om_{C/S}(-p(S)))$, and hence, by
duality, of $R^1\pi_*(\OO_C(p(S)))$. Conversely, (2) implies the similar condition for every fiber, hence, by part (i),
we get the pointwise version of (3), and the condition (3) itself follows.
\ed

We will also use the following version of the residue theorem for singular curves.

\begin{lem}\label{residue-lem} 
Let $\pi:C\to S$ be a  flat projective morphism of relative dimension $1$, with connected reduced geometric
fibers, and let
$p_1,\ldots,p_n:S\to C$ be disjoint sections such that
$\pi$ is smooth near each $p_i$. Then for any $\eta\in H^0(C\setminus\{p_1,\ldots,p_n\},\om_{C/S})$ one
has 
$$\sum_i \Res_{p_i}(\eta)=0.$$
\end{lem}

\Pf . Let $D=\sum n_ip_i$ be an effective relative divisor such that
$\eta\in H^0(C,\om_{C/S}(D))$. Note that $\und{\Ext}^i(\OO_D,\om_{C/S})=0$ for $i\neq 1$ and
$\und{\Ext}^1(\OO_D,\om_{C/S})$ is supported at $p_1(S)\cup\ldots\cup p_n(S)$.
Therefore, we have the relative Serre duality pairing
$$\tau:\pi_*\und{\Ext}^1(\OO_D,\om_{C/S})\otimes \pi_*(\OO_D)\to R^1\pi_*\om_{C/S}\simeq\OO_S.$$
By definition, its restriction $\tau_1$ to the section $1$ of $\pi_*(\OO_D)$ is obtained by applying the functor
$$R\pi_*R\und{\Hom}(\cdot,\om_{C/S})$$ to the projection $\OO_C\to\OO_D$.
Therefore, the exact sequence
$$0\to \OO_C(-D)\to \OO_C\to \OO_D\to 0$$
shows that the composition
$$\pi_*(\om_{C/S}(D))=\pi_*\und{\Hom}(\OO_C(-D),\om_{C/S})\rTo{\de} \pi_*\und{\Ext}^1(\OO_D,\om_{C/S})\rTo{\tau_1} \OO_S$$
is zero. The required formula will follow once we compute local contributions to $\tau_1(\de(\eta))=0$ at each point $p_i$,
corresponding to the canonical decomposition
$$\und{\Ext}^1(\OO_D,\om_{C/S})\simeq \bigoplus_{i=1}^n \und{\Ext}^1(\OO_{n_ip_i},\om_{C/S}).$$
Now the exact sequence
$$0\to \OO_C(-n_ip_i)\to \OO_C\to \OO_{n_ip_i}\to 0$$
induces an identification 
$$p_*\und{\Ext}^1(\OO_{n_ip_i},\om_{C/S})\simeq p_*(\om_{C/S}(n_ip_i)/\om_{C/S}),$$
compatible with $\de$,
and we claim that the projection to $\OO_S$ induced by $\tau_1$ is given by the residue at $p_i$. Indeed, 
we can replace $C$ by an open neighborhood of $p_i$, smooth over $S$, in which case this is standard.
\ed

Assume $(C,p_\bullet)$ satisfies (1) and (2), and let us fix a nonzero generator $\omega \in H^0(C, \omega_C)$. Note that by Lemma \ref{diff-restriction-lem},
$\omega$ does not vanish at any of the marked points $p_i$.
Below we will use this generator to define residues of rational functions at points $p_i$.

The cases $n=1,2$ require special attention depending on the characteristic of
$k$, which we will come back to later (see Sections \ref{case2} and \ref{case1}). 
For now, we assume that $n \geq 3$. 

Note that we have $H^1(C,\OO_C(p_i))=0$, hence for each $i\neq j$ we have $h^1(\OO_C(p_i+p_j))=0$ and $h^0(\OO_C(p_i+p_j))=2$.
By the residue theorem (see Lemma \ref{residue-lem}), for each $i\neq j$ there exist $h_{ij}\in H^0(C,\OO(p_i+p_j))$
such that
\[ \Res_{p_i}(h_{ij}\om)=1 \ \text{ and }\ \Res_{p_j}(h_{ij}\om)=-1. \] 
Let us fix a choice of $h_{1i}$ for $i=2,\ldots,n$.
% and then set \[ h_{ij}=h_{1j}-h_{1i}. \] 
For distinct $i,j\ge 2$ we set
\[ c_{ij}=h_{1i}(p_j). \]
The functions $1, h_{12}, h_{13},\ldots, h_{1n}$ form a basis of $H^0(C,\OO(p_1+\ldots+p_n))$.

Set $D=p_1+\ldots+p_n$. It is easy to see that for each $N\ge 2$ the functions
 \[ 1, h^r_{12}h_{13}, h^m_{12}, \ldots, h^m_{1n}, \ \text{ where }\  1\le m\le N, 1\le r\le N-1, \]
form a basis of $H^0(C,\OO_C(ND))$. 
Let us normalize $h_{12}$ and $h_{13}$ (adding a constant to them) by requiring that
\[ h_{12}(p_3)=h_{13}(p_2)=0, \]
i.e., $c_{23}=c_{32}=0$. 
Then $h_{12}h_{13}$ is regular near $p_2$ and $p_3$, hence, for distinct indices $i,j\ge 2$ 
the element $h_{1i}h_{1j}-h_{12}h_{13}$ belongs to $H^0(C,\OO_C(p_1+p_i+p_j))$. 
Looking at the residues at $p_i$ and $p_j$ we see that
\begin{equation}\label{hh-eq}
	h_{1i}h_{1j}-h_{12}h_{13}=c_{ij}h_{1j}+c_{ji}h_{1i}+d_{ij}
\end{equation}
for some constants $d_{ij}$.
Similarly,
\[ h_{12}h_{13}^2-h_{12}^2h_{13}\in H^0(C,\OO_C(2p_1+p_2+p_3)). \]
Hence, we should have a relation of the form
\begin{equation}\label{h-cubic-eq}
	h_{12}h_{13}^2-h_{12}^2h_{13}=ah_{12}h_{13}+bh_{12}+ch_{13}+d.
\end{equation}

Thus, setting $x_i=h_{1i}$ for $i=2,\ldots,n$, we get relations of the form
\begin{equation}\label{normal-form-1-eq}
	x_ix_j=x_2x_3+c_{ij}x_j+c_{ji}x_i+d_{ij}, \ \text{ for } i<j, (i,j)\neq (2,3),
\end{equation}
\begin{equation}\label{normal-form-1b-eq}
	x_2x_3^2=x_2^2x_3+a x_2x_3+bx_2+cx_3+d.
\end{equation}
Let us normalize the choices of $h_{1i}$ by requiring that
\begin{equation}\label{c-normalization-eq}
	c_{32}=c_{i3}=0\ \text{ for } \ i\ge 2, i\neq 3.
\end{equation}
A standard application of the Gr\"obner basis technique (where we
order the variables by $x_2 < x_3 < \ldots < x_n$ and use the
degree-lexicographical order on monomials) gives the following result.

\begin{lem}\label{Grobner-lem}
	Consider the algebra $A$ over a commutative ring $R$ with generators $x_2,\ldots,x_n$ (where $n\ge 3$)
	and the defining relations \eqref{normal-form-1-eq}, \eqref{normal-form-1b-eq}, where
	$c_{ij}, d_{ij}, a, b, c, d\in R$ and \eqref{c-normalization-eq} is satisfied. Then the elements
	\[ 1, x_2^mx_3, x_i^m, \ \text{  for } i\ge 2, m\ge 1 \]
	form an $R$-basis in $A$ if and only if the following relations hold:
	\begin{align} \begin{split}\label{Grobner-relations-eq}
			a &=c_{3i}-c_{2i}-c_{i2},\\
                        d_{ij} &=d_{2i}=-c,\\
	                d_{3i} &=b-c_{3i}c_{i2},\\
		        d &=-c_{3i}c-c_{2i}d_{3i},\\
                     	c_{ij}+c_{ji} &=c_{3i}+ c_{2j}+c_{j2},\\
			c_{ji} c_{2i}  &=c_{3i}c_{2i}+c_{j2}c_{2i}+c_{ji}c_{2j}  + d_{2j} ,\\
                        c_{ji}c_{i2}+c_{ij}c_{j2} &= c_{i2}c_{j2} + b + c,\\
			c_{ik} c_{jk} &= c_{ij}c_{jk}+c_{ji}c_{ik}+c_{3k}c_{2k}+d_{ij},		\end{split}		\end{align}
	
	where $i,j,k\ge 4$ are distinct. \ed
\end{lem}

In particular, from equations \eqref{Grobner-relations-eq}, we can see that for
$n\ge 4$ all the coefficients can be expressed in terms of the coefficients
$a$, $b$, $c$, $c_{i2}$, $c_{2i}$ and $c_{ij}$, where $4\le i<j$, which satisfy
further polynomial relations. We can rewrite these relations as follows.

\begin{defi} \label{affsch}
	For $n \geq 3$, let $U_n$ be the affine scheme over $\Z$ defined
	by the equations \eqref{c-normalization-eq} and
	\eqref{Grobner-relations-eq} on coordinates $c_{ij}, d_{ij}, a, b, c,
	d$.  \end{defi}

 \begin{prop}\label{better-equations-prop} (i) 	For $n=3$ we have $U_3\simeq \A^4_\Z$ with coordinates $a, b, c, d$. 
	
	For $n\ge 4$
	 let us set for each $i\ge 4$
	 \[ c_i=c_{2,i}, \ \ \ov{c}_i=c_{i,2}. \]
	 Let us also set $\ov{c}=b+c$.	 Then for $n\ge 4$ the ring of functions on $U_n$ is generated by 
	 $$a, c, \ov{c}, (c_i, \ov{c}_i)_{4\le i\le n}, (c_{ij})_{4\le i<j\le n},$$
         with the defining relations
	 	 \begin{equation}\label{U-n-eq}
		 \begin{array}{l}
(c_j-c_i)c_{ij}=(a+\ov{c}_i+c_j+\ov{c}_j)c_j-c, \\ (\ov{c}_i-\ov{c}_j)c_{ij}=(a+c_i+\ov{c}_i+c_j)\ov{c}_i-\ov{c} \ \text{ for } 4\le i<j, \\
			 c_{ij}c_{jk}+c_{ji}c_{ik}-c_{ik}c_{jk}+(a+c_k+\ov{c}_k)c_k=c,\\
			 \end{array}
	 \end{equation} 
	 for $4\le i<j<k$, where 
\[ c_{ji}=a+c_i+\ov{c}_i+c_j+\ov{c}_j-c_{ij} \]
	 In particular,
 \[ U_4\simeq \A^5_\Z \] with coordinates $a, c, \ov{c}, c_4, \ov{c}_4$, and
 \[ U_5\simeq \A^6_\Z \] with coordinates $a, c_4, \ov{c}_4, c_5, \ov{c}_5, c_{45}$.
	 
	 \noindent
	 (ii) Consider the morphism $U_{n+1}\to U_n$ forgetting the coordinates $c_{n+1}$, $\ov{c}_{n+1}$
	 and $c_{i,n+1}$.
	 Let also $\CC_n\to U_n$ be the affine family of curves defined by
	 \eqref{normal-form-1-eq}, \eqref{normal-form-1b-eq}.
	 Then the map
	 \begin{equation}\label{C-n-substitution-eq}
		 x_2\mapsto c_{n+1}, \ x_3\mapsto a+c_{n+1}+\ov{c}_{n+1}, x_i\mapsto c_{i,n+1}
	 \end{equation}
	 defines an isomorphism $U_{n+1}\xrightarrow{\sim} \CC_n$ over $U_n$.
 \end{prop}

 \Pf . (i) Eliminating the variables $b$, $d_{ij}$, as well as $c_{3i}$, we can rewrite
 \eqref{Grobner-relations-eq} as
 \begin{equation}
	 \begin{array}
		 {l} 
		 c_{ij}+c_{ji}=a+c_i+\ov{c}_i+c_j+\ov{c}_j,\\
		d+ac=(a+c_i+\ov{c}_i)c_i\ov{c}_i-\ov{c}c_i-c\ov{c}_i,\\
                 c=c_{ji}c_{j}+c_{ij}c_i-c_{j}c_i,\\
                \ov{c}=c_{ji}\ov{c}_i+c_{ij}\ov{c}_j-\ov{c}_{j}\ov{c}_{i}, \nonumber\\
					 c_{ij}c_{jk}+c_{ji}c_{ik}-c_{ik}c_{jk}+(a+c_k+\ov{c}_k)c_{k}=c.
	 \end{array}
 \end{equation}
  The last equation for a triple $(i,j,k)$ together with the other equations imply the similar equation for any permutation of $(i,j,k)$.
 Thus, eliminating in addition $c_{ji}$ for $j<i$ using the first set of equations above, we get the relations in \eqref{U-n-eq} from the last three sets of equations. 
 
 It is easy to see that rewriting $(\ov{c}_i-\ov{c}_j)(c_j-c_i)c_{ij}$ in two
 ways using \eqref{U-n-eq} we get the equation
 \[ (a+c_i+\ov{c}_i)c_i\ov{c}_i-\ov{c}c_i-c\ov{c}_i = (a+c_j+\ov{c}_j)c_j\ov{c}_j-\ov{c}c_j-c\ov{c}_j \]
for $i\neq j$.  Hence, we can also eliminate $d$ using the second set of equations above.  

\  

\noindent
 (ii) We have to compare the equations \eqref{normal-form-1-eq}, \eqref{normal-form-1b-eq} defining $\CC_n$ with the equations \eqref{U-n-eq} defining $U_{n+1}$. The case $n=3$ is easy so let us assume $n \geq 4$. Then, the equation \eqref{normal-form-1b-eq} follows from \eqref{normal-form-1-eq} for the pairs $(i,j)=(2,4)$ and $(3,4)$, Namely,
  the latter equations have form
  \[ x_2x_4=x_2x_3+\ldots, \ \ x_3x_4=x_2x_3+\ldots. \]
  Multiplying the first equation with $x_3$, the second---with $x_2$, and subtracting, we get exactly
  \eqref{normal-form-1b-eq}. Next, using the relations \eqref{Grobner-relations-eq} we can rewrite the equations
  \eqref{normal-form-1-eq} as 
  \begin{align}\begin{split}\label{rel-curve-eq}
  x_2x_i=x_2x_3+c_ix_i+\ov{c}_ix_2-c, & \ 4\le i\\
  x_3x_i=x_2x_3+(a+c_i+\ov{c}_i)(x_i-\ov{c}_i)+\ov{c}-c,  & \ 4\le i\\
  x_ix_j=x_2x_3+c_{ij}x_j+c_{ji}x_i-c,  & \ 4\le i<j.
  \end{split}
  \end{align}
  It remains to observe that the equations \eqref{rel-curve-eq} after the substitution \eqref{C-n-substitution-eq}
  match the equations \eqref{U-n-eq} with $j=n+1$ in the first two equations and $k=n+1$ in the third equation.
  \ed

Note that there is a natural $\G_m$-action on $U_n$ that comes from rescaling
the non-vanishing section $\omega$ of the dualizing sheaf on $C$. Correspondingly, the
degree of the variables $x_i$ are 1, the constants $c_{ij}$ have degree 1, $d_{ij}$ have degree 2, $a$ has degree 1, $b$ and $c$ have degree 2, and $d$ has degree 3. This action will play an important role in comparing this picture
to the moduli of $A_\infty$-structures in Section \ref{Ainf-sec}.

  \begin{cor}\label{tangent-cor} Assume $n\ge 5$. For any field $k$ the
	  dimension of the Zariski tangent space of $U_n\times \Spec(k)$ at the
	  origin is equal to $(n-1)(n-2)/2$.  The functions (of weight $1$ with respect to the $\G_m$-action) 
	  $$a,(c_i,\ov{c}_i)_{4\le i}, (c_{ij})_{4\le i<j}$$
	  form a minimal set of generators of the algebra of functions on $U_n$.
	  \end{cor}
  \Pf . We have that the maximal ideal $\mathfrak{m}_0$ of the local ring at the origin is generated by $a, c, \ov{c}, c_i, \ov{c}_i, c_{ij}$. We see from the defining relations that $c$ and $\ov{c}$ can be expressed in terms of the other generators. We get no additional linear dependences between the remaining generators in 
  $\mathfrak{m}_0/ \mathfrak{m}_0^2$. An easy count gives the result.
  \ed

  \begin{cor}\label{Hilbert-cor} 
  For $n\ge 5$ the graded ring $\OO(U_n)$ is free over $\Z$ with the Hilbert series
	  \[ h_{U_n}(t)=\frac{1}{(1-t)^{n+1}}\cdot \prod_{k=3}^{n-3}(1+kt+t^2).\]
  \end{cor}
  
  \Pf . Set $h_n=h_{U_n}$.
  Recall that we have a basis of $\OO(\CC_n)$ as an $\OO(U_n)$-module given by
  certain monomials in $x_i$. Since $\CC_n=U_{n+1}$, this implies that
  \[ h_{n+1}=(1+(n-1)t+nt^2+nt^3+\ldots)h_n=\frac{1+(n-2)t+t^2}{1-t}\cdot h_n. \]
  Since $\OO(U_5)$ is the ring of polynomials in $6$ variables of degree $1$, we have
  $h_5=(1-t)^{-6}$, and the assertion follows by induction on $n$.
  \ed

Sometimes the normalization \eqref{c-normalization-eq} is not convenient. An alternative is to consider for $n\ge 5$
the affine space $V_n$ (over $\Z$) with coordinates $(c_{ij})$ where $2\le i,j\le n$, $i\neq j$, subject to the linear relations
$$c_{ij}+c_{ji}-c_{ik}-c_{ki}-c_{jk}-c_{kj}+c_{kl}+c_{lk}=0$$
for any distinct $i,j,k,l$. We have a free action of $\G_a^{n-1}$ on $V_n$ such that $(a_i)_{i\ge 2}$ acts by
$$c_{ij}\mapsto c_{ij}+a_i,$$
and the quotient $\ov{V}_n$ is still an affine space. 
Choosing $h_{1i}\in H^0(C,p_1+p_i)$ with $\Res_{p_1}(h_{1i}\om)=1$, and setting $c_{ij}=h_{1i}(p_j)$ gives a well defined
point of $\ov{V}_n$ (since the only ambiguity is to add a constant to each $h_{1i}$).
Note that this construction is compatible with the action of $S_{n-1}$ permuting the points $p_2,\ldots,p_n$ and indices
$2,\ldots,n$.
Now Proposition \ref{better-equations-prop} gives the following result.

\begin{cor}\label{no-normalization-cor} 
The above construction gives a closed embedding $U_n\hra \ov{V}_n$ for each $n\ge 5$, so that
the projection $U_{n+1}\to U_n$ is compatible with the projection $\ov{V}_{n+1}\to \ov{V}_n$ omitting
$c_{ij}$ with $i=n+1$ or $j=n+1$.
\end{cor}

\begin{ex}\label{wheel-ex-1} Let us consider the wheel of $n\ge 3$ projective lines $C_1\cup \ldots\cup C_n$,
	  where $1\in C_i$ is glued to $0\in C_{i+1}$ and $p_i=\infty\in C_i$.
	  Let $u_i$ denote the natural parameter on $\A^1\sub C_i$.
	  We have a global section $\om$ of the dualizing sheaf on $C$ given by
	  \[ \om|_{C_i}=\frac{du_i}{u_i(1-u_i)}=-d\ln(1-u_i^{-1}), \ i=1,\ldots,n.\]
Then we have:
	  \[ x_2=h_{12}=\begin{cases} u_1 &\text{ on }C_1,\\ 1-u_2 &\text{ on }C_2,\\ 0 &\text{ on }C_i, \ \ i>2.\end{cases}
	  \]
	  and for $j\ge 3$,
	  \[ x_j=h_{1j}=\begin{cases} u_1-1 &\text{ on }C_1,\\ 0 &\text{ on }C_i, \ \ 1<i<j,\\ -u_j &\text{ on }C_j,\\
		  -1 &\text{ on }C_i, \ \ i>j.
	  \end{cases}
  \]
	  Now one easily checks that the defining relations of $C\setminus\{p_1,\ldots,p_n\}$ become
	  \begin{equation}\label{wheel-relations}
		  \begin{array}{l}
			  x_2x_j=x_2x_3, \ \ 3\le j,\\
			  x_ix_j=x_2x_3-x_j, \ \ 3\le i<j,\\
			  x_2x_3^2=x_2^2x_3-x_2x_3.
		  \end{array}
	  \end{equation}

  \end{ex}

\subsection{Case $n=2$, $\cha(k)\neq 2$}
\label{case2} 

In the $n=2$ case we always assume that {\bf char$(k)\neq 2$}. 

Let $t_1$ be a formal parameter at $p_1$ such that $\Res_{p_1} \om/t_1=1$.
Then there exists a function $f_1\in H^0(C,\OO_C(2p_1))$, such that at $p_1$, \[ f_1\equiv \frac{1}{t_1^2} \mod k [[t_1]] \] (we use the fact that $\Res_{p_1} f_1\om=0$).
Note that $f_1$ is unique up to adding a constant.
  Then the functions
  \[ 1, f_1h_{12}^r, h_{12}^m, \ \text{ where } 1\le m\le N, 0\le r\le N-2 \]
  form a basis of $H^0(C,\OO_C(ND))$.
  Let us also normalize $h_{12}$ so that at $p_1$
  \[ h_{12}\equiv \frac{1}{t_1} \mod t_1 k[[ t_1]]. \]
  Then $f_1-h_{12}^2$ is regular near $p_1$.
  Now let us normalize $f_1$ so that
  \[ (f_1-h_{12}^2)(p_1)=0. \]
  Then
  $f_1(f_1-h_{12}^2)\in H^0(C,\OO_C(p_1+2p_2))$, so we have
  \begin{equation}\label{fh-eq}
  f_1^2-f_1h_{12}^2=\alpha (f_1-h_{12}^2)+\beta h_{12}+\gamma
  \end{equation}
  for some constants $\alpha,\beta,\gamma$.
  
  The analog of Lemma \ref{Grobner-lem} in this case states that for any choice of $\alpha,\beta,\gamma$ in a commutative ring
  $R$ the $R$-algebra $A$ generated by $x$ and $y$ subject to the defining relation
  \begin{equation}\label{normal-form-2-eq}
	  y^2-yx^2= \alpha (y-x^2)+ \beta x+ \gamma
  \end{equation}
  has $(x^m, x^my)$ as an $R$-basis. In fact, \eqref{normal-form-2-eq} is simply the unfolding
  of the tacnode singularity
  \begin{equation}\label{tacnode-eq}
	  y^2-yx^2=0.
  \end{equation}
  
  We extend Definition \ref{affsch} by letting $U_2 = \A^3_{\Z[1/2]}$ be the affine space generated by $\alpha,\beta$, and $\gamma$. Note that we have a natural $\G_m$-action such that $\deg(x)=1$, $\deg(y)=\deg(\alpha)=2$, $\deg(\beta)=3$, $\deg(\gamma)=4$. 
  
 \begin{ex}\label{wheel-ex-2} In the case when $C$ is the wheel of $n=2$ projective lines
	  and $2$ is invertible (see Example \ref{wheel-ex-1}) we can use
	  \begin{equation}\label{h12-part-eq}
		  h_{12}=\begin{cases} u_1-\frac{1}{2} &\text{ on } C_1,\\ \frac{1}{2}-u_2 &\text{ on }C_2,\end{cases}
	  \end{equation}
	  \[ f_1=\begin{cases} (u_1-\frac{1}{2})^2 &\text{ on } C_1,\\ \frac{1}{4} &\text{ on }C_2.\end{cases} \]
	  Thus, the equation \eqref{fh-eq} in this case takes form
	  \[ f_1^2-f_1h_{12}^2=\frac{1}{4}(f_1-h_{12}^2). \]
  \end{ex}

\subsection{Case $n=1$, $\cha(k)\neq 2,3$}\label{case1}

The case $n=1$ corresponds to the classical family of Weierstrass curves. 
We assume that $\cha(k)\neq 2, 3$. 
Then we can choose a formal parameter $t$ at $p=p_1$ such that 
$$\om\equiv dt \mod t^4k[[t]]\cdot dt,$$
where $\om$ is a global section of the dualizing sheaf of $C$.
The condition $h^1(C,\OO_C(p))=0$ implies 
$h^0(C,\OO_C(mp))=m$ for $m\ge 1$. 
Hence, there exists a non-constant function (unique up to adding a constant) $x \in H^0(C, \OO_C(2p))$ such that at $p$, 
\[ x \equiv \frac{1}{t^2} \mod  k [[t]]. \] 
We can normalize it by adding a constant so that 
\[ x \equiv \frac{1}{t^2} \mod t^2  k [[t]] \]
(the coefficient of $t$ is zero since $\Res_p (x^2\om)=\Res_p (x^2dt)=0$). 
Similarly, since $h^0(C, \OO_C (3p))= 3$, there exists a unique function $y \in H^0(C,\OO_C(3p))$ such that 
\[ y \equiv \frac{1}{t^3} \mod t k[[t]]. \] 
Then $y^2 - x^3 \in H^0(C, \OO_C(2p))$, so we have
\begin{equation}\label{cusp}
y^2 - x^3 = \delta x + \epsilon
\end{equation} 
for some constants $\delta$ and $\epsilon$.
This also works in families over an affine base (see e.g., \cite[Lem.\ 1.2.1]{P-ainf}).
Again, \eqref{cusp} (with $\delta$ and $\epsilon$ viewed as independent variables) is simply the unfolding of the cusp singularity 
$$y^2-x^3=0.$$ 

We extend Definition \ref{affsch} to this case by letting $U_1 = \A^2_{\Z[1/6]}$ be the affine space generated by $\delta$ and $\epsilon$. Note that we have a natural $\G_m$-action such that $\deg(x)=2$, $\deg(y)=3$, and $\deg(\delta)=4$, $\deg(\epsilon)=6$. 

\subsection{Moduli spaces}

Let $S$ be any scheme. Let $\UU_{1,n}$ be the (non-separated) stack of flat,
proper, finitely presented morphisms $\CC \to S$ from an algebraic space $\CC$
together with $n$ sections, whose geometric fibers are reduced, connected
curves of arithmetic genus 1. It is shown by
Jack Hall in \cite[appendix B]{smyth}  that $\UU_{1,n}$ is an algebraic stack, locally of finite type over
$\Spec \Z$.  

\begin{defi} 
	For $n\geq 1$, we define the moduli stack $\UU^{sns}_{1,n}$ to be the open substack of $\UU_{1,n}$ consisting of curves $C$ of arithmetic genus $1$ with $n$ distinct smooth marked points $p_1,\ldots,p_n$ such that
\begin{enumerate} 
\item $h^0(\OO_C(p_i))=1$ for all $i$ and, 
\item $\OO_C(p_1+\ldots+p_n)$ is ample.
	\end{enumerate} 

We also denote by $\wt{\UU}^{sns}_{1,n}\to\UU^{sns}_{1,n}$ the $\G_m$-torsor associated with a choice of a generator of the one-dimensional space $H^0(C,\om_C)$, where $\om_C$ is a dualizing sheaf.

\end{defi}

Note that the condition (1) requires that $p_i$ are non-special divisors, which is an open condition.  

\begin{thm}\label{moduli-curves-thm}
	For $n\ge 3$
	the moduli stack $\wt{\UU}^{sns}_{1,n}$ is isomorphic to the affine scheme $U_n$ over $\Spec(\Z)$
	(see Def. \ref{affsch}), so that the open affine part of the universal curve $C\setminus D$,
	where $D=p_1+\ldots+p_n$, gets identified with the curve $\CC_n\to U_n$ given by
	equations \eqref{normal-form-1-eq} (or equivalently, \eqref{rel-curve-eq} for $n\ge 4$). 

	In the case $n=2$ (resp., $n=1$) the moduli stack $\wt{\UU}^{sns}_{1,n}$ over $\Spec(\Z[1/2])$ (resp., $\Spec(\Z[1/6])$)
	is isomorphic to $U_2=\A^3$ (resp., $U_1=\A^2$), 
	so that the affine universal curve $C\setminus D$ is given by \eqref{normal-form-2-eq}
	(resp., \eqref{cusp}).	
	
	These isomorphisms are compatible with the $\G_m$-actions described above.
\end{thm}

\Pf . This is similar to \cite[Thm.\ 1.2.3]{P-ainf}. We follow the line of argument given there. For simplicity let us assume that
$n \geq 3$---the cases $n=1,2$ can be analyzed similarly. Using the relative version (over an affine base $\Spec(R)$)
of the constructions that led to Lemma \ref{Grobner-lem} and Proposition \ref{better-equations-prop}, we can 
associate with a family $\pi : C \to \Spec R$ in $\wt{\UU}^{sns}_{1,n}$ an $R$-point of $U_{n}$.
Thus, we obtain a functor $\wt{\UU}^{sns}_{1,n} \to U_n$. 

Conversely, let $(a,c,\ov{c}, c_i, \ov{c_i}, c_{ij})$, where $4 \leq i
<j$, represent an $R$-point of $U_n$. We consider the
corresponding algebra $A$ over $R$ with generators $x_2, \ldots, x_n$ and
defining relations \eqref{rel-curve-eq}, or equivalently,
\eqref{normal-form-1-eq}, \eqref{normal-form-1b-eq}, where $b,d,d_{ij}$ and $c_{ij}$ for $i>j$
are determined from \eqref{Grobner-relations-eq}. 
Let $(F_m A)_{m\ge 0}$ be the increasing filtration on $A$ associated with the generators $x_2,\ldots,x_n$,
so that $F_0 A = R\cdot 1$, $F_1 A/F_0A = R x_2 \oplus R x_3 \oplus \ldots \oplus R x_n$ and for $m \geq 2$,  
$F_m A=(F_1 A)^m$. Lemma \ref{Grobner-lem} implies that for $m\ge 2$, $F_m A/F_{m-1}A$
is freely generated over $R$ by $x_2^{m-1} x_3$ and $x_i^m$ for $2 \leq i \leq n$. 
Let \[ \mathcal{R}A := \bigoplus_{m\geq 0} F_{m} A \] be the associated Rees algebra, and 
consider the corresponding projective scheme over $\Spec(R)$, 
\[ C = \Proj (\mathcal{R} A). \]

Let $T \in F_{ 1} A$ be the element corresponding to $ 1 \in F_0 A \sub F_1 A$, and let $D = (T=0)$ be the corresponding
divisor in $C$. Then one has an isomorphism $C \backslash D \cong \Spec A$ and the complementary closed set is given by \[ D \cong \Proj ( \bigoplus_{m \geq 0} F_{ m}A/ F_{m-1} A). \] 

The graded algebra $\mathcal{R}A$ is the quotient of the polynomial ring $ R[T, X_2,\ldots, X_n]$ by the homogenization of the equations \eqref{normal-form-1-eq}, \eqref{normal-form-1b-eq},
\begin{align}\label{homogen}
	X_i X_j -X_2 X_3 &=c_{ij} T X_j+c_{ji} T X_i +d_{ij} T^2, \text{\ for\ } i<j, (i,j) \neq (2,3), \\
\label{homogen2} 
	X_2 X_3^2-X_2^2 X_3 &= a T X_2 X_3 + b T^2 X_2 +c T^2 X_3 + d T^3.
\end{align}

We have $n$ sections $p_i : \Spec R \to D$ cut out by:
\begin{align*}
	p_1 &:   T=0, X_2 = X_3 = \ldots = X_n,  \\
	p_i &:   T=0, X_i \neq 0, X_2 = X_3 = \ldots = X_{i-1} = X_{i+1}= \ldots = X_n = 0 \text{ \ for\ } 2 \leq i \leq n.
\end{align*}
Equations \eqref{homogen}, \eqref{homogen2} easily imply that 
$$D=(T=0)=\sqcup_{i=1}^n \im(p_i).$$

Over a point $ s\in \Spec R$, the fibre $C_s$ is given by the
same equations over the residue field $k(s)$ of $s$. As $F_m A$ is a free
$R$-module, we can easily compute the Hilbert function of $C_{s}$ to be: \[
	h_{C_s} (m) = \dim (F_{ m} A \otimes_R k(s)) = mn \text{\ for\ } m
	\geq 1 \] Hence, $C_s$ is a reduced, connected, degree $n$ curve in
	$\mathbb{P}_{k(s)}^{n-1}$ of arithmetic genus 1. We conclude that $ \pi : C
	\to \Spec R$ has relative dimension 1. 
	%This also gives us that $D$ is a
	%disjoint union of the sections $\im (p_i)$, as we have $n$ distinct
	%sections that lie in the hyperplane section $(T=0)$.
	 Furthermore, $\OO_{\mathbb{P}^{n-1}_{R}} (1)|_C = \OO_C(1)$ is ample, so the divisor $D$ is ample.

Note that as $A$ is a free $R$-module, the morphism $\pi : C \backslash D \to \Spec R$
is flat. For each $j\ge 2$ consider the distinguished open subset $V_j = \Spec A_j\sub C$, where
$A_j$ is the degree 0 part of the localization $(\mathcal{R}A)_{X_j}$.
Since $D\sub V_2\cup\ldots\cup V_n$, it suffices to check that each $A_j$ is flat over $R$.
We know that $(A_j)_{T/X_j}$ is flat over $R$ since $\pi : C \backslash D \to \Spec R$ is flat. 
On the other hand, $A_j/ (\frac{T}{X_j}) \cong R\oplus R$ since $V_j\cap (T=0)$ 
is the disjoint union of the section $p_1$ and $p_j$, and so $A_j/(\frac{T}{X_j})$ is also flat over $R$.  Applying \cite[Lem.\ 1.2.4]{P-ainf}
 we conclude that $A_j$ is flat over $R$. 

Next, let us show that the projection  $\pi: C \to \Spec R$ is smooth near
$p_1,\ldots, p_n$. We can work with a geometric fiber of $\pi$, i.e., assume that $R$ is an algebraically closed field.
Let us show first the smoothness at $p_j$, where $j > 1$. The maximal ideal $\mathfrak{m}_{p_j}$ of the local ring at $p_j$
is generated by $T/X_j, X_i/X_j$ for
$i \neq j$. Suppose first that $j\geq 4$. Over the open set $V_j$ 
we can write using \eqref{homogen},
$$ \frac{X_i}{X_j} = \frac{X_2}{X_j} \frac{X_3}{X_j} + c_{ij} \frac{T}{X_j}  + c_{ji} \frac{T}{X_j} \frac{X_i}{X_j}  + 
d_{ij} \left(\frac{T}{X_j}\right)^2,$$
hence we have 
\[ \frac{X_i}{X_j} \equiv  c_{ij} \frac{T}{X_j} \mod \mathfrak{m}^2_{p_j}. \] 
which implies that $\mathfrak{m}_{p_j}/\mathfrak{m}^2_{p_j}$ is generated by a single element, the image
of $T/X_j$. Hence, $C$ is smooth at $p_j$. For $j=2$, the same argument works, except for the fraction $X_3/X_2$:
here we need to use that over $V_2$ one has
$$\frac{X_3}{X_2} = \left(\frac{X_3}{X_2}\right)^2 - a \frac{T}{X_2} \frac{X_3}{X_2} - b \left(\frac{T}{X_2}\right)^2 - 
c \left(\frac{T}{X_2}\right)^2 \frac{X_3}{X_2} - d \left(\frac{T}{X_2}\right)^3$$
(see \eqref{homogen2}), hence $X_3/X_2 = 0 \in \mathfrak{m}^2_{p_2}$. A similar
argument works for $j=3$.

Now let us prove smoothness at $p_1$. Note that $p_1$ lies in all of the open sets $V_i$ for
$i \geq 2$, so we can work on $V_2\cap\ldots\cap V_n$.  
The maximal ideal $\mathfrak{m}_{p_1}$ of the local ring at $p_1$ is generated by 
$T/X_2$ and $X_i/X_2-1$ for $3 \leq i\leq n$.  
For each $j\geq 4$ we have from \eqref{homogen}
$$\frac{X_j}{X_2} - \frac{X_3}{X_2} \equiv c_{2j} \frac{T}{X_2}\frac{X_j}{X_2}  + 
c_{j2} \frac{T}{X_2}\mod \mathfrak{m}_{p_1}^2\equiv (c_{2j}+c_{j2})\frac{T}{X_2}\mod \mathfrak{m}_{p_1}^2.$$ 
On the other hand, dividing \eqref{homogen2} by $X_2^2X_3$ we get 
$$\frac{X_3}{X_2}-1 \equiv a \frac{T}{X_2} \mod \mathfrak{m}_{p_1}^2.$$ 
Again, we conclude that $T/X_2$
generates $\mathfrak{m}_{p_1}/ \mathfrak{m}_{p_1}^2$.  

Next, we have to specify a choice of a global section of $\omega_C$. We know
that our family is of arithmetic genus 1. Thus, we can determine a global
1-form $\omega$ using the requirement that \[ \Res_{p_1}(x_2 \omega)=1 \]
(where $x_2 = X_2/ T$ is a function on $C$ with simple poles at $p_1$ and
$p_2$).  The same reasoning also works in a family. Therefore, we conclude that
$\pi: C \to \Spec R$ defines an object of the moduli stack
$\wt{\UU}^{sns}_{1,n}$.   

Finally, observe that $H^0(C, \OO(mD))$ can be identified with $F_{ m} A$
inside the algebra $A$ of functions on $C \backslash D$, by analyzing the polar
conditions at the marked points $p_1, \ldots, p_n$. Furthermore, we have \[ \Res_{p_1} (x_i\omega) =1  \text{\ and\ } \Res_{p_i} (x_i \omega) = -1 \text{\ for \ } 2 \leq i \leq n , \] 
as $x_i-x_2$ is regular at $p_1$ and $x_i = X_i/T$ has (simple) poles only at $p_1$ and $p_i$. Hence, it is clear
that the functors that we constructed from $U_n$ to $\wt{\UU}^{sns}_{1,n}$ and from $\wt{\UU}^{sns}_{1,n}$ to $U_n$ are inverses of each other. 
\ed

\begin{cor}\label{flat-cor} For $n\ge 3$ the scheme $\wt{\UU}^{sns}_{1,n}$ is flat over $\Z$.
\end{cor}

\Pf . This follows from the identification $\wt{\UU}^{sns}_{1,n}\simeq U_n$, Corollary \ref{Hilbert-cor} (for the case $n\ge 5$), 
and isomorphisms $U_3\simeq \A^4_\Z$, $U_4\simeq \A^5_\Z$.
\ed

\end{subsection}

\subsection{Comparison with Smyth's moduli spaces}\label{Smyth-sec}

For each $m$, $1\le m<n$, Smyth defined in \cite{Smyth-I} the notion of {\it
$m$-stability} for $n$-pointed curves of arithmetic genus $1$ and showed that
the corresponding moduli stack $\ov{\MM}_{1,n}(m)$ is an irreducible projective
Deligne-Mumford stack over $\Spec{ \Z[1/6] }$. Below we recall the definition.

Let us denote by $C_{1,n}$ the singular curve corresponding to the point in the
moduli space $\wt{\UU}^{sns}_{1,n}\simeq U_n$ where all coefficients are zero. Thus, the
curve $C_{1,n}$ has $n$ smooth points at infinity $p_1,\ldots,p_n$ such that
$C_{1,n}\setminus \{p_1,\ldots,p_n\}$ is given by the equations
$$x_ix_j=x_2 x_3$$ for $n\ge 4$ (where the indices $i<j$ vary in $[2,n]$), by
the equation $$x_2x_3^2=x_2^2x_3$$ for $n=3$, by the tacnode equation
\eqref{tacnode-eq}: $$y^2 - yx^2=0 $$ for $n=2$, and finally for $n=1$ we get the ordinary cusp
$$y^2 -x^3 = 0.$$ 

Alternatively, we can describe $C_{1,n}$ for $m\ge 3$ as the union of $n$ generic lines
passing through one point in the projective space $\P^{n-1}$. The arising
singularity at this point is called the {\it elliptic $n$-fold point} (for all $n\ge 1$), which is a
Gorenstein singularity (see \cite[Prop.\ 2.5]{Smyth-I}). We also refer to $C_{1,n}$ as the {\it elliptic $n$-fold curve}. 

Next, let us recall (see \cite[Lem.\ 3.1]{Smyth-I}) that each Gorenstein curve $C$ of arithmetic genus $1$
has the so-called {\it fundamental decomposition} 
\begin{equation}\label{fund-decomposition}
C=E\cup R_1\cup\ldots\cup R_k
\end{equation}
where $E$, called the {\it minimal elliptic subcurve} of $C$, is a connected subcurve of arithmetic genus $1$
with no disconnecting nodes, 
$R_i$ are nodal curves of arithmetic genus $0$, 
$R_i\cap R_j=\emptyset$ for $i\neq j$, and $R_i\cap E$ is a single point which is a node of $C$.

The notion of $m$-stability for an $n$-pointed curve of arithmetic genus $1$ consists of the following three conditions:

\begin{defi} \label{mstab} A curve $(C,p_1,\ldots,p_n)$ of arithmetic genus $1$ (where all marked points are smooth and distinct) is said to be {\it $m$-stable} if

(1) The curve $C$ has only nodes and elliptic $l$-fold points, $l\le m$, as singularities; \\
(2) If $E\sub C$ is the minimal elliptic subcurve then 
$$|E\cap \ov{C\setminus E}|+|E\cap\{p_1,\ldots,p_n\}|>m;$$ \\
        (3) One has $H^0(C,\TT_C(-p_1-\ldots-p_n))=0$.
\end{defi}

Note that in the original definition the condition (2) is required to hold for any
connected subcurve of arithmetic genus $1$. The fact that it is enough to require this condition for the minimal elliptic subcurve follows from \cite[Lem.\ 3.5]{Smyth-I}.

In \cite[Lem.\ 3.10]{Smyth-I}, Smyth shows that $m$-stability is
a deformation-open condition, hence one can define $\ov{\MM}_{1,n}(m)$ as an open substack
of the stack $\UU_{1,n}$ of $n$-pointed curves of arithmetic genus 1. In \cite[Th.\
3.18]{Smyth-I}, Smyth verifies that $\ov{\MM}_{1,n}(m)$ is an irreducible proper Deligne-Mumford stack over $\Spec \Z[1/6] $. 

The requirement for $6$ to be
invertible is caused by the pathology that the cuspidal curve $C_{cusp}=C_{1,1}$ has extra vector fields in
characteristics $2$ and $3$, while the tacnode $C_{tn}=C_{1,2}$ has extra vector fields in characteristic
$2$, which can lead to non-trivial infinitesimal automorphisms of $(C,p_1,\ldots,p_n)$.  
Let us explain this phenomenon in more detail.
Recall that the affine parts of $C_{cusp}$ and $C_{tn}$ are the plane curves given by
\begin{align*}
C_{cusp}\setminus \{p\} &: y^2 = x^3 \\ 
C_{tn}\setminus\{p_1,p_2\} &: y^2 = yx^2
\end{align*}
It is easy to see that in $\cha(k)\neq 2$ we have $\text{rk} H^0(C_{tn}, \TT_{C_{tn}})= 3$, while in
$\cha(k) \neq 2,3$ we have $\text{rk} H^0(C_{cusp}, \TT_{C_{cusp}})=2$ (see \cite[Prop.\ 2.3]{Smyth-I}). 
In Lemma \ref{char-2-3-lem} below we write out explicitly the extra vector fields in $\text{char}(k)=2, 3$ that
prevent condition (3) of $m$-stability to hold for these curves unless there are sufficiently many
marked points on them.

Let $C$ be either the cusp or the tacnode. Let $\nu: \tilde{C} \to C$ be the normalization map. One has $\tilde{C} = \P^1$ for the cusp curve, and it is the disjoint union of two $\P^1$'s for the tacnode. Any vector field on $C$ can be restricted to $C \backslash Sing(C) \cong \tilde{C} \backslash \nu^{-1}(Sing(C))$ and then extended to a rational vector field on $\tilde{C}$. This leads to a natural inclusion map: 
\[ \TT_C \hookrightarrow \nu_* \TT_{\tilde{C}} \otimes K(\tilde{C}) \]

\begin{lem}\label{char-2-3-lem}
	Consider the normalization map $\nu : \tilde{C} \to C_{cusp}$ defined in affine coordinates by $k[x,y]/(y^2-x^3) \to k[t^2,t^3] \subset k[t]$ sending $(x,y) \to (t^2,t^3)$. 

	If $\text{char}(k)=2$, the natural inclusion $\TT_{C_{cusp}} \to \nu_* \TT_{\tilde{C}} \otimes K(\tilde{C})$ is generated by the vector fields \[ t^2 \partial_t, t \partial_t, \partial_t, \frac{1}{t^2} \partial_t \] which are images of the sections $x^2 \partial_y, y \partial_y, x \partial_y, \partial_y$ of $\TT_{C_{cusp}}$.
	
	In particular, there is no non-zero vector field on $C_{cusp}$ which vanishes at $5$ distinct points.
	
	If $\text{char}(k)=3$, the natural inclusion $\TT_{C_{cusp}} \to \nu_* \TT_{\tilde{C}} \otimes K(\tilde{C})$ is generated by the vector fields \[ t^2 \partial_t, t \partial_t, \partial_t, \frac{1}{t} \partial_t \] which are images of the sections $x^2 \partial_x, x \partial_x, y \partial_x, \partial_x$ of $\TT_{C_{cusp}}$.
        
	In particular, there is no non-zero vector field on $\TT_{C_{cusp}}$ which vanishes at $4$ distinct points.

	Consider the normalization map $\nu : \tilde{C} \to C_{tn}$ defined in affine coordinates by $k[x,y] / (y^2-yx^2) \to k[t] \oplus k[s]$ sending $x \to (t,s)$ and $y \to (t^2, 0)$. 
	
	If $\text{char}(k)=2$, the natural inclusion $\TT_{C_{tn}} \to \nu_* \TT_{\tilde{C}} \otimes K(\tilde{C})$ is generated by vector fields \[ t^2 \partial_t, s^2 \partial_s, t \partial_t + s \partial_s, \partial_t + \partial_s \] which are images of the sections $y \partial_x, (x^2-y) \partial_x, x \partial_x, \partial_x$ of $\TT_{C_{tn}}$.        

	In particular, there is no non-zero vector field on $\TT_{C_{tn}}$ which vanishes at $5$ distinct points.
\end{lem}	
\Pf . This is an extension of \cite[Prop.\ 2.3]{Smyth-I} to the case of $\text{char}(k)=2,3$ and follows from a similar calculation as given there. 
\ed

Therefore, to include characteristics $2$ (resp., $3$) one can either throw away curves which have cusps or tacnodes with
fewer than $5$ (resp., $4$) marked points which in general may result in losing properness, or 
to relax the condition (3) which would lead to an algebraic stack which is not a Deligne-Mumford stack. In our current study, we will only be concerned with the moduli stacks $\ov{\MM}_{1,n}(n-1)$, and we choose the option
of relaxing the condition (3). Thus, we propose the following version of these stacks over $\Z$.

\begin{defi} \label{modified} The stack $\ov{\MM}_{1,n}^\infty$ is the moduli stack of
(reduced, connected projective) pointed curves $(C,p_1,\ldots,p_n)$ of arithmetic genus $1$ (where all marked points are smooth and distinct) such that

(1)' The curve $C$ has only nodes and elliptic $l$-fold points, $l<n$, as singularities; \\
(2)' $C$ has no disconnecting nodes (i.e., it coincides with its minimal elliptic subcurve); \\
(3)' Every irreducible component of $C$ contains at least one marked point.
\end{defi}

We will see in Theorem \ref{curve-moduli-thm} below that $\ov{\MM}_{1,n}^\infty$ for $n\ge 3$ (resp., $\ov{\MM}_{1,2}^\infty$)
is proper over $\Spec(\Z)$ (resp., $\Spec(\Z[1/2])$) and
that in fact it is a projective scheme for $n\ge 5$. 
The following result compares it to Smyth's moduli stack
$\ov{\MM}_{1,n}(n-1)$ over $\Spec(\Z[1/6])$.

\begin{prop} Let $(C,p_1,\ldots,p_n)$ be a reduced, connected projective curve of arithmetic genus $1$ with smooth
distinct marked points, over an algebraically closed field $k$. If $(C,p_1,\ldots,p_n)$ is $(n-1)$-stable then
conditions (1)'--(3)' are satisfied. The converse is true under one of the following additional assumptions:

(a) $\cha(k)\neq 2,3$;\\
(b) $\cha(k)\neq 2$ and $n\ge 4$;\\
(c) $n\ge 5$.
\end{prop}

\Pf . It is easy to see that for conditions (2) with $m=n-1$ and (3) to be satisfied 
$C$ has to coincide with its minimal elliptic subcurve. By \cite[Lem.\
3.3]{Smyth-I}, this implies that we have the following possibilities for the underlying curve $C$:

\noindent
(i) a smooth elliptic curve;

\noindent
(ii) an irreducible rational nodal curve;

\noindent
(iii) a wheel of $\P^1$'s;

\noindent
(iv) $C=C_{1,m}$, the elliptic $m$-fold curve.

In each of these cases one can easily see that condition (3) implies (3)'.
Conversely, if (1)'--(3)' are satisfied then to check (3) we use the above classification of minimal subcurves,
together with \cite[Cor.\ 2.4]{Smyth-I} (note that the argument of \cite[Cor.\ 2.4]{Smyth-I} works also in the case $n=2$, $\cha(k)=3$) and Lemma \ref{char-2-3-lem}.
\ed

\begin{cor} There is a natural isomorphism $\ov{\MM}_{1,n}(n-1)\simeq \ov{\MM}_{1,n}^\infty\times\Spec(\Z[1/6])$.
\end{cor}

The following observation will also be useful later.

\begin{lem}\label{CY-cor} 
For every curve $C$ corresponding to a point of $\ov{\MM}_{1,n}^\infty$, 
the dualizing sheaf $\om_C$ is isomorphic to $\OO_C$.
\end{lem}

\Pf . This follows from the fact that $C$ coincides with its minimal elliptic subcurve and from \cite[Lem.\ 3.3]{Smyth-I}.
\ed

On the other hand, we can consider the GIT stability for the action of $\G_m$ on $\wt{\UU}^{sns}_{1,n}$.
Since the degrees of all the coordinates $c_{ij},d_{ij},a,b,c,d$ on $\wt{\UU}^{sns}_{1,n}$ are positive,
the GIT-semistable points (with respect to the identity character of $\G_m$)
are exactly the points where not all coordinates vanish. Hence, this is precisely the open subscheme
$$\wt{\UU}^{sns}_{1,n}\setminus\{C_{1,n}\}\sub \wt{\UU}^{sns}_{1,n}$$ 

We denote by $\la$ the standard tautological line bundle over $\UU^{sns}_{1,n}$ with the fiber $H^0(C,\om_C)$ over
$(C,p_\bullet)$.

\begin{thm}\label{curve-moduli-thm} 
 Let us work over $\Spec(\Z)$ for $n\ge 3$, over $\Spec(\Z[1/2])$ for $n=2$, and over
$\Spec(\Z[1/6])$ for $n=1$.
One has an isomorphism of stacks
	$$\ov{\MM}_{1,n}^{\infty} \simeq
	(\wt{\UU}^{sns}_{1,n}\setminus\{C_{1,n}\})/\G_m$$ compatible with the
	inclusion into the stack $\UU_{1,n}$ of all $n$-pointed curves of
	arithmetic genus $1$.  Hence, we have isomorphisms
	\begin{align*}
	\ov{\mathcal{M}}_{1,1}^{\infty} &\simeq \P(4,6) \ \text{ over } \Spec(\Z[1/6]), \\ 
	\ov{\MM}_{1,2}^{\infty} &\simeq\P(2,3,4)  \ \text{ over } \Spec(\Z[1/2]), \\
	\ov{\MM}_{1,3}^{\infty} &\simeq\P(1,2,2,3), \\
	\ov{\MM}_{1,4}^{\infty} &\simeq\P(1,1,1,2,2),\\
	\ov{\MM}_{1,5}^{\infty} &\simeq\P^5, 
\end{align*}
where $\P(d_1,\ldots,d_k)$ denotes the weighted projective stack. 
	
	For $n\ge 6$, the stack
	$\ov{\MM}_{1,n}^\infty$ is isomorphic to the $n$-dimensional irreducible projective variety in $\P^{n(n-3)/2}_\Z$ given by
	the equations obtained from \eqref{U-n-eq} by eliminating variables $c$
	and $\ov{c}$. 
	The corresponding line bundle $\OO(1)$ is isomorphic to $\la$.
For every field $k$, the variety $\ov{\MM}_{1,n}^\infty\times \Spec(k)$
	is still irreducible. 
	\end{thm}

\begin{rem} The identifications of $\ov{\MM}_{1,n}^\infty$ for $n\le 4$ show that it is not
a Deligne-Mumford stack over $\Spec(\Z)$ (non-\'etale automorphism groups occur in characteristics $2$ and $3$). 
However, it is a ``tame stack'' in the sense of \cite{AOV}, which
seems to be a better notion in positive and mixed characteristics.
\end{rem}

\Pf . First, we are going to check that any point $[(C,p_1,\ldots,p_n)] \in
\wt{\UU}^{sns}_{1,n}$, different from $C_{1,n}$, satisfies conditions (1)'--(3)'. Note that (3)' holds by definition.

To see that $C$ is Gorenstein, we need to check that the dualizing sheaf
$\omega_C$ is locally free. As was mentioned above, we know that $C_{1,n}$ is
Gorenstein by \cite[Prop.\ 2.5]{Smyth-I}, i.e. $\omega_{C_{1,n}}$ is
locally free. It follows that any curve $C$ in a neighborhood of $C_{1,n}$ is
Gorenstein. Finally, note that the $\G_m$-action brings an arbitrary point in
$\wt{\UU}^{sns}_{1,n}$ to an isomorphic curve in a neighborhood of $C_{1,n}$. 

Next, using the fact that $h^0(C,\OO(p_i))=1$, we check condition (2)'. Indeed, each
subcurve $R_i$ in the fundamental decomposition \eqref{fund-decomposition}
should have at least one marked point $p_j$ since
$\OO(p_1+\ldots+p_n)$ is ample. But then for such a point we necessarily have
$h^0(C,\OO(p_j))\ge h^0(R_i,\OO(p_j))=2$, since the arithmetic genus of $R_i$
is zero. Hence, $C$ coincides with its minimal elliptic subcurve.

It follows that 
$C$ can be either a smooth elliptic curve, an irreducible
rational nodal curve, a wheel of $\P^1$'s, or the elliptic
$m$-fold curve with $m < n$ (see \cite[Lem.\ 3.3]{Smyth-I}). This immediately gives (1)'.

Conversely, if conditions (1)'--(3)' are satisfied then we can easily check that $h^0(C,\OO(p_i))=1$ for each marked point $p_i$.
Thus, we see that $\ov{\MM}_{1,n}^\infty$ is an open substack in $\UU^{sns}_{1,n}$ obtained by throwing away
the elliptic $n$-fold point $C_{1,n}$.

The explicit identification of $\ov{\MM}_{1,n}^{\infty}$ for $n\leq 5$ follows from the identification of $U_n$ given in Proposition \ref{better-equations-prop} for $n=3,4,5$ and $n=1,2$ follows from the identification given by \eqref{cusp} and \eqref{normal-form-2-eq} (see Theorem \ref{moduli-curves-thm}). 

The irreducibility of $\ov{\MM}_{1,n}^\infty$ (resp., $\ov{\MM}_{1,n}^{\infty}\times \Spec(k)$) 
for $n\ge 6$ can be proved by the same method as in 
\cite{Smyth-I} (using the fact that the relevant curves are smoothable).

The identification $\OO(1)\simeq\la$ follows from the fact that the pull-backs of both line bundles to
$\wt{\UU}^{sns}_{1,n}\setminus\{C_{1,n}\}$ have natural trivializations, and the corresponding $\G_m$-actions on the trivial 
line bundle are both given by the identity character $\G_m\to\G_m$.
\ed

In Proposition \ref{Gr-prop} below we will also prove an isomorphism
$$\ov{\MM}_{1,6}^{\infty}\simeq Gr(2,5),$$
where $Gr(2,5)$ denotes the Grassmannian of 2-dimensional subspaces of a vector
space of dimension 5.

\begin{cor}\label{Sn-action-cor} For $n\ge 5$ the natural action of $S_n$ on $\ov{\MM}_{1,n}^{\infty}$ is 
induced by some automorphisms of the projective space $\P^{n(n-3)/2}$.
\end{cor}

\Pf . Indeed, this follows from the fact that the projective embedding $\ov{\MM}_{1,n}^{\infty}\hra\P^{n(n-3)/2}$ is given by
the $S_n$-equivariant line bundle $\la$.
\ed

\begin{rem} In terms of the isomorphism $\ov{\MM}_{1,n}^\infty\simeq (U_n\setminus\{0\})/\G_m$ (which is obtained
by combining Theorems \ref{moduli-curves-thm} and \ref{curve-moduli-thm}), the action of 
generators of $S_n$ can be easily described. Namely, the subgroup of permutations fixing $1,2,3$ acts
by natural permutations of the coordinates $c_i,\ov{c}_i, c_{ij}$ (see Proposition \ref{better-equations-prop}).
The transpositions $(13)$, $(23)$ and $(24)$ act by the following involutions: 
\begin{align*} 
(13):&\ a\mapsto a, \ c\leftrightarrow \ov{c}, \ c_i\leftrightarrow \ov{c}_i, \ c_{ij}\mapsto c_{ji}=a+c_i+\ov{c}_i+c_j+\ov{c}_j-c_{ij},\\
(23):&\ a\mapsto -a,\ c\mapsto c-\ov{c}, \ \ov{c}\mapsto-\ov{c}, \ c_i\mapsto a+c_i+\ov{c}_i, \ \ov{c}_i\mapsto -\ov{c}_i, \
c_{ij}\mapsto c_{ij}-\ov{c}_i,\\
(24):&\ a\mapsto-a-2c_4-2\ov{c}_4, \ c\mapsto\ov{c}-\ov{c}_4(a+c_4+\ov{c}_4), \ \ov{c}\mapsto c-c_4(a+c_4+\ov{c}_4), \\
&\ c_i\leftrightarrow c_{4i}, \ \ov{c}_i\leftrightarrow c_{i4}, \ c_{ij}\mapsto c_{ij}, \
\text{ where } i,j\ge 5.
\end{align*}
These formulas can be checked using the equations of the relative curve 
\eqref{rel-curve-eq}. A more transparent way to see the action of the subgroup $S_{n-1}$ fixing $1$, is via
the identification of the affine space containing $U_n$ with the quotient $\ov{V}_n$ considered in
Corollary \ref{no-normalization-cor}.
  \end{rem}   
  
We can now reprove the result of Smyth that
the moduli stack $\ov{\MM}_{1,n}(m)$ is smooth if and only if $m\le 5$, and get some additional information
on the singularities of $\ov{\MM}_{1,n}(m)$ for $m>5$. 

First, we need a bit of deformation theory.
Let us fix a field $k$, and for each $m\ge 1$ let us consider the deformation functor $\Def_m$ of the $m$-fold elliptic singularity over $k$ (defined on local Artinian algebras with the residue field $k$).
We denote by $B_m$ the base of the formal miniversal deformation of the elliptic $m$-fold singularity.
Thus, $B_m$ is the formal spectrum of $R_m$, a complete local Noetherian ring with the residue field $k$, and we have a 
formally \'etale morphism of functors 
$$h_{R_m}\to \Def_m$$
where $h_{R_m}$ is the representable functor corresponding to $R_m$. 

Let $(C,p_\bullet,v_\bullet)$ be a point of $U_n=\wt{\UU}_{1,n}^{sns}$, and let $\hat{\OO}$ 
be the completion of its local ring. Let also $\hat{U}_n$ be the formal completion of $U_n$ at this point,
i.e., the formal spectrum of $\hat{\OO}$.
Let $q_1,\ldots,q_r$ be all the non-nodal singular points of $C$, where $q_i$ is the elliptic $m_i$-fold singularity.
By \cite[Lem.\ 2.1]{Smyth-II}, the natural projection
\begin{equation}\label{def-sing-pts-map}
p:h_{\hat{\OO}}\to \prod_{i=1}^r \Def_{m_i}
\end{equation}
is formally smooth (we use the fact that the deformations of nodal singularities are unobstructed). 

\begin{lem}\label{def-Bm-lem} 
The morphism \eqref{def-sing-pts-map} 
factors through a formally smooth morphism 
$\wt{p}:h_{\hat{\OO}}\to h_R$, where
$$v:h_R\to \prod_{i=1}^r \Def_{m_i},$$
is the miniversal deformation with $R=\hat{\ot} R_{m_i}$.
Thus, $\wt{p}$ corresponds to a formally smooth morphism
$$\hat{U_n}\to \prod_{i=1}^r B_{m_i}.$$
%One has an isomorphism of $R$-algebras
%$\hat{\OO}\simeq R[[t_1,\ldots,t_s]]$.
\end{lem}

\Pf .
By definition, the morphism $v$ is formally \'etale. In particular, it is formally smooth, so we can lift (non-uniquely) the projection
$p$ to a local homomorphism $R\to \hat{\OO}$, 
%(with the trivial residue field extension)
so that we have a commutative triangle 
\begin{equation}
\begin{diagram}
h_{\hat{\OO}} \\
\dTo{\wt{p}} &\rdTo{p}\\
h_R&\rTo{v}& \prod_{i=1}^r \Def_{m_i}
\end{diagram}
\end{equation}
It then follows that the morphism $\wt{p}$ is formally smooth.
% i.e., $\hat{\OO}$ is a formally smooth $R$-algebra. 
Indeed, since $v$ is \'etale, and $p$ is smooth, the morphism $\wt{p}$ induces a surjection 
on tangent spaces, hence, we can apply \cite[Lem.\ 4.5.3]{P-ainf}.
%Since the residue fields of $R$ and $\hat{\OO}$ are the same, 
%the last assertion follows from \cite[Prop.\ 1.5]{SGA1-III}???
\ed

The analog of Lemma \ref{def-Bm-lem} also holds with $U_n$ replaced by $\ov{\MM}_{1,n}(m)$. This
implies that $\ov{\MM}_{1,n}(m)$ is smooth if and only if
$B_{m'}$ is smooth for $m'\le m$. Since we have a formally smooth
morphism from the completion of $U_m$ at zero to $B_m$, this is equivalent to smoothness of $U_{m'}$ at zero
for $m'\le m$.  Now we recall that $U_m$ is smooth for $m\le 5$ (see Proposition \ref{better-equations-prop}).
On the other hand, $U_6$ is not smooth at zero, since it
is $7$-dimensional but the dimension of the tangent space is $10$ (see
Corollary \ref{tangent-cor}). Thus, we recover the fact that $\ov{\MM}_{1,n}(m)$ is smooth if and only if $m\le 5$
(see \cite[Cor.\ 4.17]{Smyth-II}).

In a similar vein we have the following results. In the rest of this subsection we work over an algebraically closed field.

\begin{prop}\label{codim-6-smooth-prop} 
The moduli stacks $\ov{\MM}_{1,n}^{\infty}$ and the schemes $U_n$ are smooth in codimension $\le 6$.
\end{prop}

\Pf . Use the stratification by the singularity type (see \cite[Cor.\ 2.4]{Smyth-II}) and the fact that $B_m$ is smooth for $m\le 5$.
\ed

\begin{prop}\label{Gor-norm-prop} 
For $n\ge 3$ the scheme $U_n\simeq\wt{\UU}^{sns}_{1,n}$ is Gorenstein and
regular in codimension $1$, hence, normal.
\end{prop}

\Pf . Regularity in codimension $1$ follows from Proposition \ref{codim-6-smooth-prop}. To prove the Gorenstein property
we use the induction on $n$. We have seen in Proposition \ref{better-equations-prop} that
$U_n$ is smooth for $n\le 5$. Assume the assertion is true for $U_n$. By Theorem \ref{moduli-curves-thm}
and by Proposition \ref{better-equations-prop}(ii), we know that
$U_{n+1}$ is open in the universal curve over $U_n=\wt{\UU}^{sns}_{1,n}$. 
So the morphism $U_{n+1}\to U_n$ is flat, with Gorenstein fibers (see the proof of Theorem \ref{curve-moduli-thm})
and Gorenstein base. Hence, $U_{n+1}$ is Gorenstein too (see \cite[Ch.\ V, Prop.\ 9.6]{hartshorne}). 
\ed

\begin{cor}\label{normal-cor} 
For $n\ge 5$ the projective scheme $\ov{\MM}_{1,n}^{\infty}$ is normal and arithmetically Gorenstein, hence
Gorenstein. 
\end{cor}

\Pf . Indeed, the homogeneous coordinate ring of $\ov{\MM}_{1,n}^{\infty}$ is exactly the ring of functions on $U_n$.
\ed

\begin{cor}\label{normal-cor2} 
The stacks $\ov{\MM}_{1,n}(m)$ are normal and Gorenstein for $1\le m<n$. 
\end{cor}

\Pf . Recall that the Gorenstein property of a local Noetherian ring can be checked after passing to its completion 
(see \cite[Thm.18.3]{Matsumura-CRT}).
Also, if $R\to S$ is a local formally smooth homomorphism of complete Noetherian rings, then $R$ is Gorenstein if
and only if $S$ is Gorenstein (this follows from the results of \cite{AF}).
Thus, by Lemma \ref{def-Bm-lem}, to check the Gorenstein property
it suffices to check that the base $B_m$ of formal miniversal deformation of
the elliptic $m$-fold point is Gorenstein. But we have a formally smooth morphism from the completion of $U_m$ at zero
to $B_m$ (by the same Lemma), and $U_m$ is Gorenstein by Proposition \ref{Gor-norm-prop}.
Hence, $B_m$ is also Gorenstein. On the other hand, as in Proposition \ref{codim-6-smooth-prop},
we see that stacks $\ov{\MM}_{1,n}(m)$ are regular in codimension $1$. Hence, we conclude that they are normal.
\ed

\subsection{Rational map from $\ov{\MM}_{1,n}^{\infty}$ to $\ov{\MM}_{1,n-1}^{\infty}$ and its applications}

For each $1\le i<j\le n$ let us consider the point 
$P_{ij}=P_{ij}(n)\in\ov{\MM}_{1,n}^{\infty}$ corresponding to the elliptic $(n-1)$-fold curve
with the marked points $p_i$ and $p_j$ on the same component and exactly one marked point on each other component.

Let 
$$\pi_n:\ov{\MM}_{1,n}^{\infty}\dashrightarrow \ov{\MM}_{1,n-1}^{\infty}$$
be the rational map corresponding to omitting the marked point $p_n$.
Note that it is compatible via the isomorphism of Theorem \ref{curve-moduli-thm}
with the linear projection $U_n\to U_{n-1}$ omitting
the coordinates $c_n$, $\ov{c}_n$ and $c_{in}$, $i=4,\ldots,n-1$ (here and below we use the coordinates introduced
in Proposition \ref{better-equations-prop}(i)).
Therefore, in terms of the natural embeddings $\ov{\MM}_{1,n}^{\infty}\sub\P^{n(n-3)/2}$,
$\ov{\MM}_{1,n-1}^{\infty}\sub\P^{(n-1)(n-4)/2}$, 
the rational map $\pi_n$ is induced by the linear projection
\begin{equation}\label{lin-projection-eq}
\P^{n(n-3)/2}\dashrightarrow \P^{(n-1)(n-4)/2}
\end{equation}
along the $n-3$-dimensional
projective subspace $\P(K_n) \sub\P^{n(n-3)/2}$ where the linear subspace $K_n$ is
given by the equations $c_{ij}=0$ for $4 \leq i<j\le n-1$, $c_i=\ov{c}_i=a=0$ for $4\leq i\le n-1$.

\begin{prop}\label{blow-up-prop}
Assume that $n\ge 6$. 

\noindent
(i) The intersection $\ov{\MM}_{1,n}^{\infty}\cap \P(K_n)$ is transversal and consists of the points $P_{in}=P_{in}(n)$, $i=1,\ldots,n-1$.
More precisely, the homogeneous coordinates $c_n$, $\ov{c}_n$, $c_{in}$ ($4\le i<n$) at these points are:
\begin{equation}
\begin{array}{l}\nonumber
P_{1n}: \ \ov{c}_n=0, c_{in}=c_n\neq 0;\\
P_{2n}: \ \ov{c}_n=-c_n\neq 0, c_{in}=0;\\
P_{3n}: \ \ov{c}_n\neq 0, c_n=c_{in}=0;\\
P_{in}, i\ge 4: \ c_n=\ov{c}_n=0, c_{in}\neq 0, c_{jn}=0 \text{ for } j\neq i.
\end{array}
\end{equation}

\noindent
(i') In terms of coordinates $(c_{ij})$, where $2\le i,j\le n$, $i\neq j$ (see Corollary \ref{no-normalization-cor}),
the point $P_{1i_0}$, where $2\le i_0\le n$ is determined by the equations
\begin{align}
\begin{split}
c_{ij}=c_{ik},\\
c_{ii_0}-c_{ik}=c_{ji_0}-c_{jk},
\end{split}
\end{align}
where the indices $i,j,k$ are distinct and different from $i_0$.

\noindent
(ii) The rational map $\pi_n$ is resolved by a diagram of regular maps
\begin{equation}\label{blow-up-diagram}
\begin{diagram}
&&\ov{\CC}_{n-1}\\
&\ldTo{r_n}&&\rdTo{q_{n-1}}\\
\ov{\MM}_{1,n}^{\infty}&&&& \ov{\MM}_{1,n-1}^{\infty}
\end{diagram}
\end{equation}
where $q_{n-1}:\ov{\CC}_{n-1}\to \ov{\MM}_{1,n-1}^{\infty}$ is the universal curve.
Via the projection $r_n$, $\ov{\CC}_{n-1}$ gets identified with the blow-up of $\ov{\MM}_{1,n}^{\infty}$ at
the $n-1$ points $(P_{in})$, $i=1,\ldots,n-1$. The image of the canonical section
$\si_i:\ov{\MM}_{1,n-1}^{\infty}\to\ov{\CC}_{n-1}$ (where $i=1,\ldots,n-1$) coincides with the exceptional divisor 
$E_i=r_n^{-1}(P_ {in})$.

\noindent
(iii) One has $\pi_n(P_{ij}(n))=P_{ij}(n-1)$ for $j\le n-1$. Furthermore, $P_{ij}(n)$ is the singular point of the elliptic
$(n-2)$-fold curve $\pi_n^{-1}(P_{ij}(n-1))$. Let $C_{ij}\sub \pi_n^{-1}(P_{ij}(n-1))$ be the component containing
$p_i$ and $p_j$, and let $C_k\sub \pi_n^{-1}(P_{ij}(n-1))$ be the component containing $p_k$, where $k\neq i,j$.
Then $r_n(C_k)$ is the line connecting $P_{kn}(n)$ and $P_{ij}(n)$, while $r_n(C_{ij})$ is a conic containing
$P_{ij}(n)$, $P_{in}(n)$ and $P_{jn}(n)$.
\end{prop}

\Pf . (i) Proposition \ref{better-equations-prop}(ii) shows that with respect to the homogeneous coordinates
$y_2,\ldots,y_{n-1}$ on the projective space $\P(K_n)$ given by
$$y_2=c_n, \ y_3=c_n+\ov{c}_n, \ y_i=c_{in} \text{ for } 4\le i<n,$$
the scheme-theoretic intersection $\ov{\MM}_{1,n}^{\infty}\cap \P(K_n)\sub \P(K_n)$ is defined 
by the equations $y_iy_j=y_2y_3$, for any $2\le i<j<n$. 
This immediately implies that this intersection is transversal and consists of the following $n$ points:
\begin{equation}
\begin{array}{l}\nonumber
P'_{1n}: \ y_2=\ldots=y_{n-1}\neq 0;\\
P'_{in}, 2\le i<n: \ y_i\neq 0, y_j=0 \text{ for } j\neq i.
\end{array}
\end{equation}
It remains to show that $P'_{in}=P_{in}$. For this we use the equations \eqref{rel-curve-eq} of the affine curve associated
with each of the points $P'_{in}$. It is easy to check that this curve is a union of $n-2$ lines
and a conic $C_{in}$, which has as two points at infinity the marked points $p_i$ and $p_n$. More precisely, with respect to
the coordinates $x_2,\ldots,x_n$ on the affine part of the corresponding curve this conic component is given by
\begin{equation}
\begin{array}{l}\nonumber
C_{1n}: \ x_2=\ldots=x_{n-1}, x_2x_n=x_2^2+c_nx_n;\\
C_{2n}: \ x_3=\ldots=x_{n-1}=0, x_2x_n=c_n(x_n-x_2);\\
C_{3n}: \ x_2=x_4=\ldots=x_{n-1}=0, x_3x_n=\ov{c}_n(x_n-\ov{c}_n);\\
C_{in}, 4\le i<n: \ x_ix_n=c_{in}(x_n-x_i), x_j=0 \text{ for } j\neq i,n.
\end{array}
\end{equation}

\noindent
(i') Applying a transposition swapping $i_0$ and $n$ we can assume $i_0=n$. 
Now the result follows from (i).

\noindent
(ii) Let $r:B\to \ov{\MM}_{1,n}^{\infty}$ be the blow-up at $\ov{\MM}_{1,n}^{\infty}\cap \P(K_n)=\{P_{in} \ |\ i=1,\ldots,n-1\}$.
We can realize both $B$ and $\ov{\CC}_{n-1}$ as (reduced) subschemes of the projective bundle 
$$\P(\OO^{n-2}\oplus \OO(-1))\to \ov{\MM}_{1,n-1}^{\infty}.$$
Indeed, by (i), $B$ is a closed subscheme of the blow-up $\wt{B}$
of $\P^{n(n-3)/2}$ along the projective subspace $\P(K_n)$. The linear projection \eqref{lin-projection-eq}
extends to a regular map $\wt{B}\to\P^{(n-1)(n-4)/2}$, which can be identified with the projective bundle
$\P(\OO^{n-2}\oplus \OO(-1))$ over $\P^{(n-1)(n-4)/2}$.
On the other hand, the embedding of the relative curve $\ov{\CC}_{n-1}$ into the same projective bundle
corresponds to the surjection
$$(x_2,\ldots,x_{n-1},1): q_{n-1}^*(\OO^{n-2}\oplus \la)\to q_{n-1}^*\la(p_1+\ldots+p_{n-1}),$$
where $x_i=h_{1i}$ are rational functions defined in Sec.\ \ref{normal-forms-sec}. Here, instead
of normalizing $h_{1i}$ by their residue at $p_i$, we view them as canonical morphisms
$$h_{1i}:\la^{-1}\simeq \OO(p_i)|_{p_i}\to (q_{n-1})_*\OO(p_1+p_i-p_3), \text{ for } i\neq 3,$$
$$h_{13}:\la^{-1}\simeq \OO(p_3)|_{p_3}\to (q_{n-1})_*\OO(p_1+p_3-p_2)$$
(the isomorphisms with $\la^{-1}$ follow from Lemma \ref{diff-restriction-lem}).

Let $H\sub\P(\OO^{n-2}\oplus\OO(-1))$ be the relative hyperplane at infinity, i.e., the image of
the embedding $\P(K_n)\times \ov{\MM}_{1,n-1}^{\infty} \hra \P(\OO^{n-2}\oplus\OO(-1))$ corresponding to the
embedding of bundles $\OO^{n-2}\sub \OO^{n-2}\oplus\OO(-1)$.
The intersection $H\cap\ov{\CC}_{n-1}$ is exactly the union of the canonical sections $\si_1,\ldots,\si_{n-1}$.
Note that $\ov{\CC}_{n-1}\setminus H$ can be identified with the quotient of the affine family of curves
$$\CC_{n-1}\setminus K_n\to U_{n-1}\setminus 0$$
by the action of $\G_m$. By Proposition \ref{better-equations-prop}(ii), we deduce the equality
$$\ov{\CC}_{n-1}\setminus H=r^{-1}(\ov{\MM}_{1,n}^{\infty} \setminus \P(K_n)).$$
Passing to closures we get the equality of the subschemes in $\P(\OO^{n-2}\oplus \OO(-1))$,
$$\ov{\CC}_{n-1}=B.$$

The images of the canonical sections $\si_i$ in $H\cap \ov{\CC}_{n-1}\sub H=\P(K_n)\times \ov{\MM}_{1,n-1}^{\infty}$ are given by the equations
\begin{equation}
\begin{array}{l}\nonumber
\si_1: \ y_2=\ldots=y_{n-1};\\
\si_i, 2\le i<n: \ y_j=0 \text{ for } j\neq i,
\end{array}
\end{equation}
where $(y_i)$ are the homogeneous coordinates on $\P(K_n)$ defined in part (i).
Since the restriction of the projection $r_n:\ov{\CC}_{n-1}\to \ov{\MM}_{1,n}^{\infty} \sub \P^{n(n-3)/2}$ to $H\cap\ov{\CC}_{n-1}$ is given by
the natural projection $\P(K_n)\times \ov{\MM}_{1,n-1}^{\infty} \to\P(K_n)\sub \P^{n(n-3)/2}$, we deduce that $r_n(\si_i)=P'_{in}=P_{in}$ as claimed.

\noindent
(iii) Without loss of generality we can assume that $i=1$ and $j=n-1$ (see Corollary \ref{Sn-action-cor}). 
Then the fact that $\pi_n(P_{1,n-1}(n))=P_{1,n-1}(n-1)$ 
follows immediately from (i') and from Corollary \ref{no-normalization-cor}. 
Next, as in (i) we see that the elliptic $(n-2)$-fold curve corresponding to $P_{1,n-1}(n-1)$ is given by the equations
\begin{align*}
&y_iy_j=y_2y_3 \text{ for } 2\le i<j<n-1,\\
&y_2y_3+c_{n-1}y_{n-1}=y_iy_{n-1} \text{ for } 2\le i<n-1.
\end{align*}
Using this one can easily check the remaining assertions.
\ed

\begin{cor} For $n\ge 6$, the tangent cone to $\ov{\MM}_{1,n}^{\infty}$ at the point $P_{ij}$ is isomorphic to the affine cone over
	$\ov{\MM}_{1,n-1}^{\infty}\sub \P^{(n-1)(n-4)/2}$.
\end{cor}

\Pf . For $j=n$ this follows from Proposition \ref{blow-up-prop}(ii). The general case follows using the action of
the symmetric group $S_n$ on $\ov{\MM}_{1,n}^{\infty}$.
\ed

We now show how the diagram \eqref{blow-up-diagram} can be used in studying the geometric properties of 
our moduli spaces.

\begin{prop}\label{can-prop} Assume $n\ge 5$.

\noindent
(i) The canonical line bundle $K$ on $\ov{\MM}_{1,n}^{\infty}$ is
$K\simeq \OO(n-11)$.

\noindent
(ii) Over an algebraically closed field of characteristic zero or over $\Spec(\Z)$ one has 
\[ \Pic \ov{\MM}_{1,n}^{\infty} =\Z, \] 
and this group is generated by the class $\OO(1)=\la$.

\noindent
(iii) Let us work over an algebraically closed field of characteristic zero.
Then the variety $\ov{\MM}_{1,n}^{\infty}$ has rational singularities if and only if $n\le 11$.
\end{prop}

\Pf . For brevity we denote the maps in diagram \eqref{blow-up-diagram} as $q=q_{n-1}$, $r=r_n$.

\noindent
(i) Let $E\sub\ov{\CC}_{n-1}$ be the exceptional divisor of the blow-up $r$ (consisting of $n-1$ components).
Since $q$ is induced by the linear projection, we have
\begin{equation}\label{divisors-eq}
q^*\OO(1)\simeq r^*\OO(1)(-E).
\end{equation}
Next, we note that by Lemma \ref{CY-cor}, the relative dualizing sheaf on the universal curve satisfies
$$\om_q\simeq q^*\la=q^*\OO(1).$$
Hence, we can prove our assertion by induction in $n$.
For $n=5$ this is true since $\ov{\MM}_{1,5}^{\infty}\simeq\P^5$.
Assuming that the canonical bundle on $\ov{\MM}_{1,n-1}^{\infty}$ is $\OO(n-12)$ we get
$$K_{\ov{\CC}_{n-1}}\simeq \om_q\ot q^*\OO(n-12)\simeq q^*\OO(n-11).$$
Therefore, by \eqref{divisors-eq}, the canonical bundle of
$$\ov{\MM}_{1,n}^{\infty}\setminus\{P_{1n},\ldots,P_{n-1,n}\}\simeq \ov{\CC}_{n-1}\setminus E$$
is still $\OO(n-11)$, and the induction step follows (recall that $\ov{\MM}_{1,n}^\infty$ is normal by Corollary \ref{normal-cor}).

\noindent
(ii) First, let us work over $\C$.
We use the fact that the rational Picard group of $\ov{\MM}_{1,n}^{\infty}$ is $\Q$ (\cite[Prop. 3.2]{Smyth-II}).
Let us show now that the group $\Pic(\ov{\MM}_{1,n}^\infty)$ has no torsion. Suppose $\xi$ is a torsion line bundle
on $\ov{\MM}_{1,n}^{\infty}$. Then $L=r^*\xi$ is a torsion line bundle on $\ov{\CC}_{n-1}$ with the property
$L|_{E_i}\simeq \OO$ for $i=1,\ldots,n-1$, where $E_i$ are the components of the exceptional divisor.
Recall that by Proposition \ref{blow-up-prop}(ii), $E_i$ is the image of the $i$th canonical section of $q$.
Next, we claim that $M:=q_*(L(E_1))$ is a line bundle on $\ov{\MM}_{1,n-1}^{\infty}$. Indeed, for every fiber $C$ of $q$
the line bundle $L|_C$ is torsion, hence, it has degree zero on every irreducible component of $C$. Let $p_1=E_1\cap C$.
It is enough to show that $H^1(C,L|_C(p_1))=0$ (since then $H^0(C,L|_C(p_1))$ will be one-dimensional by Riemann-Roch).
By Lemma \ref{CY-cor}, we have $\om_C\simeq\OO_C$, so by Serre duality we need to show the vanishing of $H^0(C,L^{-1}|_C(-p_1))$. But the line bundle $L^{-1}|_C(-p_1)$
has degree $-1$ on one of the irreducible components of $C$ and degree $0$ on the remaining components, so it has
no global sections. Thus, $M=q_*(L(E_1))$ is a line bundle. Let $D$ be the effective divisor given as the vanishing locus 
of the natural map $q^*M\to L(E_1)$, so that 
\begin{equation}\label{line-bundle-section-eq}
L(E_1)\simeq q^*M(D).
\end{equation}
Then on each fiber $C$ of $q$ we have
$L|_C(p_1)\simeq \OO(D)|_C$. We claim that the unique global section of $L|_C(p_1)$ vanishes at exactly one smooth point.
This is easy to see when $C$ is either irreducible or a wheel of projective lines. Suppose now that $C$ is an elliptic $m$-fold
curve, and let $C_1\sub C$ be a component containing $p_1$. It suffices to show that we cannot have a global section of
$L|_C(p_1)$ vanishing at all the other components of $C$. Indeed, restricting such a global section to
a neighborhood of the singular point $q$, we will get a germ $f\in \OO_{C,q}$ which restricts to zero on all branches
but one 
and will have a nonzero derivative at $q$ on the remaining branch. This contradicts the explicit description of $\OO_{C,q}$
(see \cite[Sec.\ 2]{Smyth-I} and Lemma \ref{tangent-action-lem}(i) below), 
so no such global section exists. Hence, $D$ defines a section of $q$, so it is isomorphic
to $\ov{\MM}_{1,n-1}^{\infty}$. Assume first that $D$ is different from all the divisors $E_i$. 
Then the conditions $L|_{E_i}\simeq \OO$ for $i=1,\ldots,n-1$ give isomorphisms
\begin{align*}
&\OO(-1)\simeq M(q(D\cap E_1)),\\
&\OO\simeq M(q(D\cap E_i)) \ \text{ for }i\ge 2
\end{align*}
on $\ov{\MM}_{1,n-1}^{\infty}$. In particular, the Cartier divisors $q(D\cap E_2)$ and $q(D\cap E_3)$ are linearly
equivalent. Since $\Pic(\ov{\MM}_{1,n-1}^{\infty})_{\Q}=\Q$, every nonzero effective Cartier
divisor on $\ov{\MM}_{1,n-1}^{\infty}$ is ample. Since the divisors $q(D\cap E_2)$ and $q(D\cap E_3)$ do not intersect,
we deduce that $D\cap E_2=0$, so $M\simeq \OO$. But then we should have $\OO(-1)\simeq \OO(q(D\cap E_1))$ which
is a contradiction. Suppose next that $D=E_i$ with $i\ge 2$. Then restricting \eqref{line-bundle-section-eq} to
$E_1$ and to $E_i$ we get $\OO(-1)\simeq M$ and $\OO\simeq M\ot\OO(-1)$, which is a contradiction.
It follows that $D=E_1$, i.e., $L=q^*M$. Restricting to $E_1$ we get $M\simeq \OO$, hence, $L$ is trivial.

Finally, to check that $\OO(1)$ is a generator of $\Pic(\ov{\MM}_{1,n}^{\infty})$, we observe that $\ov{\MM}_{1,n}^{\infty}$
contains a projective line (with respect to the embedding given by $\OO(1)$). Indeed, we can use of the lines contained
in $\pi_n^{-1}(P_{ij}(n-1))$ (see Proposition \ref{blow-up-prop}(iii)). Hence, $\OO(1)$ is not divisible in the Picard group
of $\ov{\MM}_{1,n}^{\infty}$.

The same result over $\Spec(\Z)$ follows by the standard method (see \cite[p.\ 103]{Mumford}), using the irreducibility
of the fibers of $\ov{\MM}_{1,n}^\infty\to\Spec(\Z)$.

\noindent
(iii) First, let us check that $\ov{\MM}_{1,n}^{\infty}$ has rational singularities for $n\le 11$, by induction on $n$.
For $n\le 6$ this is true since our moduli space is smooth. Assume $\ov{\MM}_{1,n-1}^{\infty}$ has rational
singularities. Let $P$ be a point of $\ov{\MM}_{1,n}^{\infty}$ which is different from all $(P_{ij})$,
so that the corresponding curve $(C,p_1,\ldots,p_n)$ has at most 
elliptic $(n-2)$-fold points as singularities. Applying the $S_n$-action to $P$ we can assume that there is more than one marked point on the component of $C$ containing $p_n$.
Then viewing $P$ as a point on $\ov{\CC}_{n-1}$
we see that the projection $\ov{\CC}_{n-1}\to \ov{\MM}_{1,n-1}^{\infty}$ is smooth near $P$.
Thus, $P$ has an open neighborhood $U$ which is smooth over $\ov{\MM}_{1,n-1}^{\infty}$, hence
$U$ has rational singularities. It remains to check that $\ov{\MM}_{1,n}^{\infty}$ has rational singularities near
each of the points $P_{ij}$. Using the $S_n$-action, it is enough to consider the points $P_{in}$.
Note that the projection $q:\ov{\CC}_{n-1}\to \ov{\MM}_{1,n-1}^{\infty}$ is smooth near $E$, hence
$\ov{\CC}_{n-1}$ has rational singularities near $E$. Now by Lemma \ref{rat-sing-lem} below, applied to the
blow-up morphism $r$, $\ov{\MM}_{1,n}^{\infty}$ has rational singularities near $P_{in}$ if and only if
$$r_*K_{\ov{\CC}_{n-1}}\simeq K_n,$$
where $K_n$ is the canonical bundle of $\ov{\MM}_{1,n-1}^{\infty}$. As we have seen in part (i),
$$K_{\ov{\CC}_{n-1}}\simeq q^*\OO(n-11)\simeq r^*K_n((11-n)E),$$
where we used  \eqref{divisors-eq}. Thus, 
$$r_*K_{\ov{\CC}_{n-1}}\simeq K_n\ot r_*(\OO((11-n)E)).$$
Since $H^0(E,\OO_E(-i))=0$ for $i>0$, we see that
$r_*(\OO(iE))=\OO$ for $i>0$, which finishes the induction step.
The same argument for $n=12$ shows that $\ov{\MM}_{1,12}^{\infty}$ does not have rational singularities
at the special points $P_{ij}$ (since $r_*(\OO(-E))\neq \OO$). 
For any $n>12$ we can find a point $P\in \ov{\MM}_{1,n}^{\infty}$ projecting to a special point
in $\ov{\MM}_{1,12}^{\infty}$, such that an open neighborhood of $P$ projects smoothly to
$\ov{\MM}_{1,12}^{\infty}$. Hence, $\ov{\MM}_{1,n}^{\infty}$ will not have a rational singularity near $P$.
\ed

We have used the following result, well-known to the experts.
% (see e.g., \cite[Lem.\ 1, Thm.\ 3]{Kovacs} for similar statements).

\begin{lem}\label{rat-sing-lem} Let $\phi:Y\to X$ be a proper birational morphism of quasi-projective varieties over 
an algebraically closed field of characteristic zero. Assume that $Y$ has rational singularities. Then $X$ has
rational singularities if and only if $X$ is Cohen-Macaulay and $\phi_*\om_Y\simeq \om_X$.
\end{lem}

\Pf . Let $\pi:Z\to Y$ be the resolution of singularities, and let $p=\phi\circ\pi$. 
It is well known (see e.g., \cite[Lem.\ 1, Thm.\ 3]{Kovacs}) that $Y$ has rational singularities if and only if
$Y$ is Cohen-Macaulay and $\pi_*\om_Z\simeq\om_Y$.
Similarly, $X$ has rational singularities if and only if 
$X$ is Cohen-Macaulay and $p_*\om_Z\simeq \om_X$. Since $\pi_*\om_Z\simeq\om_Y$, the latter condition is
equivalent to $\phi_*\om_Y\simeq\om_X$.
\ed

\begin{rem}
Rationally (and over an algebraically closed field of characteristic zero), 
the identification of the canonical class of $\ov{\MM}_{1,n}^{\infty}$
 with $(n-11)\la$  follows from the description of 
$\ov{\MM}_{1,n}^{\infty}$ as an explicit birational contraction of $\ov{M}_{1,n}$ in \cite{Smyth-II}. 
\end{rem}

\subsection{Curves which are linear sections of $Gr(2,5)$}

In the case $n=6$ we can identify the blow-up picture of Proposition \ref{blow-up-prop} with a generic linear projection
of the Grassmannian $Gr(2,5)$ to $\P^5$.

\begin{prop}\label{Gr-prop} 
	There is an isomorphism $\ov{\MM}_{1,6}^{\infty} \cong Gr(2,5) \subset \P^9$, so that the map
	$q_5$ from the blow up of $\ov{\MM}_{1,6}^{\infty}$ at $5$ points to $\ov{\MM}_{1,5}^{\infty}\simeq\P^5$ 
	(see \eqref{blow-up-diagram}) gets identified with a generic linear
projection of $Gr(2,5)$ (that can be defined over $\Z$).
\end{prop} 

First, we consider the well known family of (degenerating) elliptic curves obtained by linear sections of $Gr(2,5)$
(see e.g., \cite[Sec.\ 2]{Coskun-Prend}). Namely, we consider a fixed linear subspace $L=\P^3\sub\P^9$, intersecting
$Gr(2,5)$ at $5$ distinct points $p_1,\ldots,p_5$. Then for every $4$-dimensional subspace $P\sub \P^9$,
containing $L$, the intersection $Gr(2,5)\cap P$ is a curve of arithmetic genus $1$, containing the points $p_1,\ldots,p_5$.

\begin{lem}\label{Gr-sec-lem} 
Let $k$ be an algebraically closed field, and consider the above picture over $k$.
Then for every $4$-dimensional subspace $P\sub\P^9$ containing $L$, the curve $C=Gr(2,5)\cap P$ with $5$ marked points
$p_1,\ldots,p_5$ is $4$-stable, i.e., defines a point of $\ov{\MM}_{1,5}^{\infty}$.
\end{lem}

\Pf . The curve $C=Gr(2,5)\cap P$ is Gorenstein and is of arithmetic genus $1$. Also, the line bundle
$\OO_C(p_1+\ldots+p_5)=\OO_C(1)$ is ample. 
By \cite[Prop.\ 2.3]{WY}, the points $p_1,\ldots,p_5$ are in general linear position.
Hence, 
$$h^0(\OO_C(p_i))=h^0(\OO_C(1)(-p_1-\ldots-\hat{p_i}-\ldots-p_5))=1$$
for $i=1,\ldots,p$. Therefore, $(C,p_1,\ldots,p_5)$ corresponds to a point of $\UU^{sns}_{1,5}$. 
By Theorem \ref{curve-moduli-thm}, it remains to check that we cannot have $C\simeq C_{1,5}$. 
Since $C$ has degree $5$, if this were the case then $C$ would have to be a union of $5$ lines $\ell_i$, connecting a point
$q\in Gr(2,5)$ with each of the five points $p_i=[U_i]$. Let $T\sub\P^9$ be the tangent space to $Gr(2,5)$ at $q$.
Then we have $\ell_i\in T$ for $i=1,\ldots,5$, hence, $P\sub T$. It is well known that the intersection $Gr(2,5)\cap T$
is $4$-dimensional (in fact, it is the cone over the Segre embedding of $\P^1\times\P^2$, see 
\cite[Sec.\ 6.1.3]{AG}). Since $P$ is a linear subspace of codimension $2$ in $T$, this implies that
$Gr(2,5)\cap P$ has dimension $\ge 2$, which is a contradiction.
\ed

\noindent 
{\it Proof of Proposition \ref{Gr-prop}.}
It is enough to choose the embedding $L=\P^3\sub\P^9$ defined over $\Z$ such that $Gr(2,5)\cap L$ consists of
$5$ distinct points (over any algebraically closed field).
Indeed, let $X\to Gr(2,5)$ be the blow up along $p_1,\ldots,p_5$. Then the fibers of the regular map $\pi:X\to\P^5$
(defined as a linear projection with the center $L$) are exactly linear sections of $Gr(2,5)$ by $4$-dimensional
subspaces containing $L$. The exceptional divisors $E_1,\ldots,E_5\sub X$ give five non-intersecting sections of $\pi$.
Thus, by Lemma \ref{Gr-sec-lem}, we get a family of $4$-stable curves, which is the pull-back
of the universal family with respect to some regular map $f:\P^5\to\ov{\MM}_{1,5}^{\infty}$. Recall that $\ov{\MM}_{1,5}^{\infty}\simeq\P^5$ so that
$\la=\OO(1)$. We claim that $f^*\OO(1)\simeq \OO(1)$. Indeed, we have to calculate the line bundle $\la$ associated with
the family $\pi:X\to\P^5$. By Lemma \ref{diff-restriction-lem}(ii), we have an isomorphism
$$\la\simeq \om_\pi|_{E_1}\simeq \OO_X(-E_1)|_{E_i}$$
for any $i=1,\ldots,5$.
But $\OO_X(-E_i)|_{E_i}\simeq\OO(1)$ since $E_i$ is the exceptional divisor of a blow-up. 
Hence, $f$ is an isomorphism and the map $\pi:X\to\P^5$ is a universal
curve over $\ov{\MM}_{1,5}^{\infty}$. Now the identification of $\ov{\MM}_{1,6}^{\infty}$ follows from Proposition \ref{blow-up-prop}.

It remains to give a subspace $\P^3\sub\P^9$ with the required properties. 
%We do this by a direct computation as follows. 
We use the homogeneous coordinates 
$$[ z_{12} : z_{13} : z_{14} : z_{15} : z_{23} : z_{24} : z_{25} : z_{34} : z_{35} : z_{45}] $$
on $\P^9$ associated with the standard basis of $\bigwedge^2\Z^5$.	
We choose the following four points on $Gr(2,5)$:
\begin{align*}
p_1 &= [1:0:-1: \ \ 0:0:1:0:0:0:0], \ \
&p_2 = [0:1:\ \ 0:-1:0:0:0:0:1:0], \\
p_3 &= [0:0:\ \ 0:\ \ 0:1:0:0:0:0:0], \ \
&p_4 = [0:0:\ \ 0:\ \ 0:0:0:0:0:0:1], 
\end{align*}
and let $L=\P^3$ be the linear subspace spanned by these four points in $\P^9$. Note that $L$ is defined by the equations
\[ z_{12} = z_{24}\ , \ z_{13} = z_{35}\ , \ z_{14} + z_{24} = 0\  ,\  z_{15} + z_{35} =0\ ,\  x_{25} =0\ ,\ x_{34} =0. \]
In particular, $L$ is defined over $\Z$.
Since $Gr(2,5) \subset \P^9$ has degree 5, one expects that $L$ intersects $Gr(2,5)$ in yet another point. Using the Pl\"ucker relations, it is easy to verify that indeed $L$ intersects $Gr(2,5)$ in exactly one other point, namely:
$$ p_5 = [1:-1:-1:1:1:1:0:0:-1:-1] .$$ \ed

\begin{rem} Over an algebraically closed field of characteristic $0$ the isomorphism of $\ov{\MM}_{1,6}^{\infty}$ with $Gr(2,5)$
	follows immediately from the fact that $\ov{\MM}_{1,6}^{\infty}$ is a smooth projective variety of dimension $6$ and degree $5$
in $\P^9$ (not contained in any hyperplane). Indeed, this is a part of the classification result of Fujita in \cite{Fujita}. 
\end{rem}

\begin{defi} Let $C$ be a projective curve of arithmetic genus $1$.
We say that $C$ is {\it $4$-prestable} if either $C$ is irreducible with at most one nodal singularity, or, $C$ is a wheel of 
$\le 5$ projective lines, or $C$ is the elliptic $m$-fold curve with $m\le 4$.
\end{defi}

Note that over an algebraically closed field, a curve $C$ is $4$-prestable if and only if there exists $5$ smooth points
$p_1,\ldots,p_5$ such that $(C,p_1,\ldots,p_5)$ is $4$-stable.

\begin{cor}\label{linear-section-cor} 
Let $k$ be an algebraically closed field. A curve $C\sub\P^4$ can be obtained as a 
linear section of $Gr(2,5)$ (in its Pl\"ucker embedding) if and only if $C$ is $4$-prestable.
\end{cor}

\begin{rem}
The fact that every $4$-stable curve $(C,p_1,\ldots,p_5)$ can be realized as a linear section of $Gr(2,5)$ can also be proved directly using the fact that such curves are arithmetically Gorenstein in the projective embedding associated with
$\OO(p_1+\ldots+p_5)$. Namely, one can mimic the construction from \cite[Sec.\ 4]{Fisher} of the Pfaffian presentation
of elliptic normal curves of degree $5$ based on Buchsbaum-Eisenbud theorem on Gorenstein ideals of codimension $3$
(see \cite{BE}). Let $I\sub S=k[x_0,\ldots,x_4]$ be the homogeneous ideal corresponding to $C$.
Then $S/I$ has a minimal free resolution of the form
$$0\to S(-5)\to S(-3)^5\rTo{f} S(-2)^5\to S\to S/I\to 0,$$
where $f$ is a skew-symmetric matrix of linear forms and $I$ is generated by the principal $4\times 4$-Pfaffians of $f$.
Interpreting $f$ as a map $\wt{f}:k^5\to \we^2(k^5)$ we get a required linear section of the Grassmannian.
In the case when $C$ is the elliptic $5$-fold curve, i.e., the union of $5$ generic lines in $\P^4$, the minimal resolution
still has the same form. However, in this case $\wt{f}$ has a one-dimensional kernel, 
so it does not give an embedding of $\P^4$ into $\P^9$. Instead, in this case the Pfaffian presentation realizes $C$ as a cone over $5$ points in $\P^3$.
\end{rem}

\section{Moduli of $A_\infty$-structures}\label{Ainf-sec}

\subsection{$A_\infty$-structures associated with curves}\label{curves-ainf-sec}

\ 

\quad For each curve $(C,p_1,\ldots,p_n)$ corresponding to a point of $\UU^{sns}_{1,n}$ 
we can consider the associative algebra $\Ext^*(G,G)$, where 
$$G=\OO_C\oplus\OO_{p_1}\oplus\ldots\oplus\OO_{p_n}.$$
It turns out that the fact that $h^0(\mathcal{O}_C (p_i)) = 1$ for all $i$
implies that the algebra $\Ext^*(G,G)$ is independent of the curve $(C,p_1,\ldots,p_n)$, up to an isomorphism
(the requirement that $\OO_C(p_1+\ldots+p_n)$ is ample is not important at the moment).

More precisely, let us denote by $A_i\in\Hom(\OO_C,\OO_{p_i})\sub\Ext^*(G,G)$ the natural generators.
Also, a choice of nonzero tangent vectors at the marked points gives canonical generators $B_i$ 
of the one-dimensional spaces $\Ext^1(\OO_{p_i},\OO_C)\sub\Ext^*(G,G)$, so that the algebra
$\Ext^*(G,G)$ is generated by $A_i$ and $B_i$
over the subalgebra $k^{n+1}$ (generated by the projectors to the summands of $G$).
The algebra structure on $\Ext^*(G,G)$ is given by the following easy computation, similar to
the one in \cite[Sec.\ 1.1]{FisPol}.

\begin{lem} Let $Q=Q_n$ be the quiver as in Figure \ref{fig1}. We identify $A_i$ with the arrow from the central vertex to the vertex $i$ and $B_i$ with the arrow in the opposite direction. A choice of nonzero tangent vectors at all $p_i$'s gives
rise to a canonical isomorphism of $k$-algebras
$$\Ext^*(G,G)=E_{1,n}:=k[Q]/J_1,$$
where $J_1$ is the ideal generated by the relations
$$B_iA_i=B_jA_j,\ \ A_iB_j=0 \ \text{ for } i\neq j.$$ \ed 
\end{lem}  

We will consider this as an isomorphism of graded algebras by declaring the gradings as $|A_i|=0$ and $|B_i|=1$. Note that this is different from the path-length grading. 

We use the convention that the paths are composed from the right.\footnote{This is opposite to the convention adopted in \cite{FisPol} and \cite{P-ainf}, where we reverse the direction of the arrows in the quiver $Q_n$.} 
For example, $B_i A_i$ is the image under the product:
\[ \Ext^*(\OO_{p_i}, \OO_C) \otimes \Ext^*(\OO_{C}, \OO_{p_i}) \to \Ext^*(\OO_{C}, \OO_C ) \] 
of the generators $B_i \in \Ext^*(\OO_{p_i}, \OO_C)$ and $A_i \in \Ext^*(\OO_C,\OO_{p_i})$.

\begin{figure}[!h] \centering
	\includegraphics[scale=1]{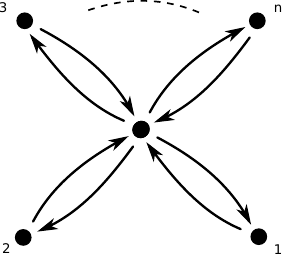}  \caption{The quiver $Q_n$}
	\label{fig1} \end{figure}

We refer to \cite[Ch. 1]{seidel} for the basics of the theory of $A_\infty$-algebras.
Some of the more relevant aspects for us are discussed in \cite{P-ainf}.

Let $R$ be a commutative ring. For a fixed graded associative $R$-algebra $A$
one can consider {\it minimal $A_\infty$-algebras} $(\mathcal{A}, \mu^* )$ over $R$
(minimality means that $\mu^1=0$) extending the given double product on $A$.
When considering $A_\infty$-algebra structures on $E_{1,n}$ we in addition assume that they are {\it unital} with
respect to the natural idempotents in $E_{1,n}$ corresponding to the vertices of the quiver $Q_n$, i.e., any operation
$\mu^l$, $l>2$, that has one of these idempotents as one of the arguments, is required to vanish. 

Working over a field $k$, it turns out that if the first Hochschild cohomology
group $HH^1(A)_{<0}$ vanishes then the set of minimal $A_\infty$-structures on the
associative algebra $A$ up to a gauge equivalence can be represented by an affine
scheme (see \cite[Cor. 4.2.5]{P-ainf}). We denote this affine scheme by
$\MM_\infty(A)$. Note that, in general, this scheme is constructed
as the inverse limit 
$$ \MM_\infty(A) = \varprojlim_d \MM_d(A) $$
of affine schemes $\MM_d(A)$ of finite type
which represents the set of minimal $A_d$-structures over $A$, 
hence is not necessarily of finite type.

As in \cite[Section 3.1]{P-ainf}  (see also, \cite[Section 5.1]{LP} for the special case $n=1$), we can use a natural dg-resolution of $\Ext^*(G,G)$ to construct an $A_\infty$-structure on $E_{1,n} \ot R$ associated with a family of curves 
$(C,p_1,\ldots,p_n)$ over $\Spec R$ in $\wt{\UU}^{sns}_{1,n}(R)$ (where $R$ is a commutative ring). This is constructed by considering endomorphisms of $G$ in a $dg$-enhancement of $D^b(C)$ and then applying the homological perturbation lemma (see \cite{Merk}, \cite{KS})
to get a minimal $A_\infty$-algebra, defined canonically up to a gauge equivalence.

Thus, we have a morphism of functors 
\begin{equation}\label{genus-1-a-inf-map}
	\wt{\UU}^{sns}_{1,n}\to \MM_\infty(E_{1,n}),
\end{equation}
where $\MM_\infty(E_{1,n})$ is the functor associating with $R$ the set of
gauge equivalence classes of minimal $A_\infty$-structures on $E_{1,n}\ot R$.
We will show later using the results of \cite[Sec.\ 4]{P-ainf}
that if we work over a field $k$ then $\MM_\infty(E_{1,n})$ is representable by an affine $k$-scheme
of finite type. Our main theorem in this section is that the map (\ref{genus-1-a-inf-map}) is an isomorphism of affine schemes. 

\begin{lem}\label{Gm-action-compatibility-lem}
(i) The map \eqref{genus-1-a-inf-map} is compatible
with the natural $\G_m$-actions, where the action on $\MM_\infty(E_{1,n})$ is given by 
\begin{equation}\label{m-n-rescaling}
(\mu^n)\mapsto (\la^{n-2}\mu^n),
\end{equation}
which is also the transformation induced by the rescaling
$B_i\mapsto \la B_i$.

\noindent
(ii) The $A_\infty$-structure on $E_{1,n}$ associated with the elliptic $n$-fold curve $C_{1,n}$ is trivial (up to a gauge
equivalence).
\end{lem}

\Pf . (i) As in \cite[Prop.\ 3.3.2]{P-ainf} (see also \cite{LP} Lemma 5.2 for $n=1$), 
one can check that it is possible to choose the homotopy operator needed to run
the homological perturbation for the dg-algebra associated with the universal curve over the affine scheme
$\wt{\UU}^{sns}_{1,n}\simeq U_n$, in a $\G_m$-equivariant way (where the $\G_m$-action is extended naturally
to the universal curve). Then, as in \cite[Prop.\ 3.3.2]{P-ainf}, one checks that the resulting higher products $\mu^n$
will have weight $n-2$ with respect to the $\G_m$-action. It is easy to check that 
the rescaling $B_i\mapsto \la B_i$ produces the same transformation \eqref{m-n-rescaling}.

\noindent
(ii) (See \cite[Prop.\ 4.4.1]{P-ainf} for a similar argument.) First, we observe that there is a natural $\G_m$-action 
on $C_{1,n}$ (induced by the $\G_m$-action on the universal curve over $U_n$), 
therefore, we get an induced $\G_m$-action on the $\Ext$-algebra
$E_{1,n}$. It is easy to see that this action is given by $\la^{deg}$ where $deg$
is the cohomological grading on $E_{1,n}$. Using the $\G_m$-equivariant
homological perturbation as in (i), we obtain a $\G_m$-equivariant $A_\infty$-structure on $E_{1,n}$.
But $\mu^n$ lowers the cohomological degree by $n-2$, so only $\mu^2$ can be nonzero.
\ed

Similarly to \cite[Prop.\ 4.4.1]{P-ainf} we get the following result connecting the associative algebra $E_{1,n}$
to the elliptic $n$-fold curve $C_{1,n}$.

\begin{prop} \label{perf} We have an equivalence of perfect derived categories
$$\Per(C_{1,n})\simeq \Per(E_{1,n})$$
inducing a $\G_m$-equivariant isomorphism 
\begin{equation}\label{genus-1-HH-isom-eq}
HH^*(C_{1,n})\simeq HH^*(E_{1,n}),
\end{equation}
so that the second grading on these spaces is given by the weights of the $\G_m$-action.
\end{prop}

\subsection{Comparison of the moduli spaces}

We have the following analog of \cite[Lem.\ 4.4.2]{P-ainf}. 
%(See also \cite[Prop 2.3]{Smyth-I}.)

\begin{lem}\label{tangent-action-lem} Let $C=C_{1,n}$, where $n\ge 2$. We work over an arbitrary field $k$, except that in the case $n=2$ we assume that $\cha(k) \neq 2$. 
	Recall that $D=p_1+\ldots+p_n$.

\noindent
(i) The algebra of functions on $C\setminus D$ can be identified with the subalgebra in
$\prod_{i=1}^n k[x_i]$ consisting of $(f_i(x_i))$ such that $f_i(0)=f_j(0)$ for $i\neq j$ and 
$$f'_1(0)=f'_2(0)+\ldots+f'_n(0).$$

\noindent
(ii) The one-dimensional space $H^1(C,\OO)$ has weight $1$ with respect to the $\G_m$-action.

\noindent
(iii) The space $H^0(C,\TT)$, where $\TT$ is the tangent sheaf, decomposes as a direct sum
$$H^0(C,\TT)=H^0(C,\TT(-D))\oplus V,$$
where $V$ is an $n$-dimensional subspace of weight $1$ with respect to the $\G_m$-action, such that
the composition
$$V\to H^0(C,\TT)\to H^0(D,\TT|_D)$$
is an isomorphism. Furthermore,
$$H^0(C,\TT(-D))=H^0(C,\TT)^{\G_m}$$ 
and this space is spanned by the derivation coming from the $\G_m$-action on
$C$. The natural map $H^0(C,\TT(nD))\to H^0(C,\TT(nD)|_D)$ is surjective for $n\ge 0$.
Also, one has 
$$H^0(C,\TT(-2D))=0.$$ 

\noindent
(iv) One has $H^1(C,\TT)=0$.
\end{lem}

\Pf . Let $C_i\simeq \P^1$, $i=1,\ldots,n$, be the irreducible components of $C$, and let $q\in C$ 
be the singular point. For $n\ge 3$ we number components in such a way that $x_i$ is a coordinate on 
$C_i\setminus\{p_i\}$ for $i\ge 2$ and $x_j=0$ on $C_i$ for $j\ge 2$, $j\neq i$. The remaining
component $C_1$ has the affine part $C_1\setminus\{p_1\}$ given by $x_i=x_j$ for $i\neq j$, and
we denote by $x_1$ the restriction of any of $x_i$'s to $C_1\setminus\{p_1\}$. 
For $n=2$, we let $x_i$ be the natural coordinates on $C_i\setminus\{p_i\}$, $i=1,2$, obtained by
restricting $x$ (so that $y=0$ on $C_1$ and $y=x^2$ on $C_2$). 
Let us also set $U=C\setminus D$,
$U_i=C_i\setminus\{p_i\}$ and $V=C\setminus\{q\}$.

\noindent
(i) Assume first that $n\ge 3$. Then 
the algebra $\OO(C\setminus D)$ has the basis $1, x_i^m, x_2^mx_3$, where $i\ge 2$, $m\ge 1$.
The projections $\OO(C\setminus D)\to k[x_i]$ for $i\ge 2$ are given by $x_i\mapsto x_i$, $x_j\mapsto 0$ for 
$j\neq i$. The projection $\OO(C\setminus D)\to k[x_1]$ sends all $x_i$ to $x_1$.
Now the assertion follows immediately by considering the images of the basis vectors.
In the case $n=2$ we have the basis $x^m, yx^m$ (where $m\ge 0$) on $\OO(C\setminus D)$.
The map
$$\OO(C\setminus D)\to k[x_1]\oplus k[x_2]$$
sends $x^m$ to $(x_1^m,x_2^m)$ and $yx^m$ to $(0,x_2^{m+2})$, and the assertion follows.

\noindent
(ii) We can compute $H^1(C,\OO)$ using the covering of $C$ by two affine open sets: $C=U\cup V$.
Thus, this group is identified with the cokernel of the map
$$H^0(U,\OO)\oplus H^0(V,\OO)\to H^0(U\cap V,\OO)=\prod_{i=1}^n k[x_i,x_i^{-1}].$$
Functions on $V=\sqcup_i C_i\setminus\{q\}$ map to collections $(P_i(x_i^{-1}))$, where $P_i$ are arbitrary
polynomials. A collection $(x_iQ_i(x_i))$, where $Q_i$ are polynomials, 
comes from an element of $H^0(U,\OO)$ if and only if $Q_1(0)=0$. Hence, the classes in $H^1(C,\OO)$
are represented by elements of the form $(ax_1)$ with $a\in k$. It remains to observe
that $(\la^{-1})^*x_1=\la x_1$.

\noindent
(iii) Let us first study derivations of the algebra $\OO(U)$. Every such derivation restricts to
a derivation of $\OO(U\setminus q)=\prod_{i=1}^n \OO(U_i)$. Assume first that $n\ge 3$.
Then we get a collection of vector fields
$v_i\in k[x_i,x_i^{-1}]\del_{x_i}$. It is easy to see that such a collection extends to a derivation of $\OO(U)$
if and only if there exists a constant $a\in k$ such that $v_i\in (ax_i+x_i^2k[x_i])\del_{x_i}$. On the other hand,
$v_i$ is regular at infinity (i.e., extends to $C_i\setminus\{q\}$) if and only if $v_i\in x_i^2k[x_i^{-1}]\del_{x_i}$.
Thus, an element of $H^0(C,\TT)$ corresponds to a collection of the form $(v_i=(ax_i+b_ix_i^2)\del_{x_i})$;
the subspace $H^0(C,\TT(-D))$ consists of $(v_i=ax_i\del_{x_i})$. On the other hand, for $n\ge 1$ the space
$H^0(C,\TT(nD))$ consists of $v_i=P_i(x_i)\del_{x_i}$, where $P_i$ are polynomials of degree $n+2$ with
$P_i(0)=0$, $P_i'(0)=a$ (independent of $i$). This easily implies all our assertions.

\noindent
(iv) As in (ii), we can use the covering $C=U\cup V$ to compute $H^1(C,\TT)$.
Thus, we just need to see that every derivation of $\OO(U\setminus q)$ is a sum of a derivation that
is regular at infinity and a derivation that extends to $U$. But this follows easily from the explicit form of
such derivations in (iii).
\ed

As in \cite[Lem.\ 4.4.3]{P-ainf}, we deduce the following results about the Hochschild cohomology of $C_{1,n}$
(for the case $n=1$, see \cite[Sec.\ 4.1, 4.2]{LP}).

\begin{cor}\label{HH-C-1n-cor} 
For $n\geq 3$, one has
$$HH^1(C_{1,n})_{<0}=0,$$
and the natural map
$$HH^2(C_{1,n})\to HH^2(U)$$
is an isomorphism, where $U=C_{1,n}\setminus D$. \\
The same conclusions hold for $n=2$ and $\cha(k) \neq 2$ (resp., for $n=1$ and $\cha(k)\neq 2,3$). 
\end{cor}

Using the isomorphism \eqref{genus-1-HH-isom-eq} we deduce the vanishing
\begin{equation}\label{HH1-vanishing-eq}
HH^1(E_{1,n})_{<0}=0,
\end{equation}
which implies the following result.

\begin{lem}\label{M-infty-repr-lem}
	For $n\geq3$ the functor $\MM_\infty(E_{1,n})$  of minimal $A_\infty$-structures on $E_{1,n}$
 is represented by an affine $k$-scheme. The same conclusion holds for $n=2$ and $\cha(k) \neq 2$
 (resp., $n=1$ and $\cha(k)\neq 2,3$). 
\end{lem}

\Pf . This follows from  \cite[Cor.\ 4.2.5]{P-ainf}) using the vanishing \eqref{HH1-vanishing-eq}.
\ed

For every $k$-scheme $X$ we denote by $\bL_X$ the cotangent complex of $X$ over $k$.
We have the following analog of Lemma 4.4.5 of \cite{P-ainf}.

\begin{lem}\label{Ext-U-lem}
For $C=C_{1,n}$ the natural maps
$$\Ext^1_C(\bL_C,\OO(-D-p_i))\to \Ext^1_C(\bL_C,\OO(-D))\to \Ext^1_U(\bL_U,\OO_U) \ \text{ and} $$
$$\Ext^2_C(\bL_C,\OO(-D-p_i))\to \Ext^2_U(\bL_U,\OO_U)$$
are isomorphisms,
where $p_i$ is any of the standard marked points on $C=C_{1,n}$.
\end{lem}
 
\Pf . 
%The proof is similar to that of \cite[Lem.\ 4.4.5]{P-ainf}.
Since $\bL_C$ is a perfect complex, it is enough to show that the maps
$$\Ext^i_C(\bL_C,\OO(-D-p_i))\to \Ext^i_C(\bL_C,\OO(-D)) \ \text{ and}$$
$$\Ext^i_C(\bL_C,\OO(nD))\to \Ext^i_C(\bL_C,\OO((n+1)D)), \ n\ge -1,$$
are isomorphisms for $i=1,2$. Using the exact sequences
$$0\to \OO(-D-p_i)\to \OO(-D)\to \OO_{p_i}\to 0,$$
$$0\to \OO(nD)\to \OO((n+1)D)\to \OO_D\to 0,$$
this reduces to the surjectivity statement in Lemma \ref{tangent-action-lem}(iii) together 
with the surjectivity of the map
$$H^0(C,\TT(-D))\to H^0(C,\TT(-D)|_{p_i})$$
which is checked similarly.
\ed

Next, we are going to compare the deformation theories of 
$\wt{\UU}^{sns}_{1,n}$ and $\MM_\infty(E_{1,n})$. This is analogous to \cite[Sec.\ 4.5]{P-ainf}, so we will be brief. 
A slight difference of our case from the one considered in \cite{P-ainf} is in
the identification of the tangent space to $\wt{\UU}^{sns}_{1,n}$ at $C=C_{1,n}$.
By Lemma \ref{diff-restriction-lem}, a choice of a nonzero global $1$-form is equivalent to a choice of a nonzero tangent
vector to one of the marked points. Thus, we can identify the tangent space to $\wt{\UU}^{sns}_{1,n}$ with 
$\Ext^1_C(\bL_C,\OO(-D-p_i))$ for any $i$---these spaces are canonically
isomorphic. In fact, by Lemma \ref{Ext-U-lem}, these spaces are naturally isomorphic to $\Ext^1_U(\bL_U,\OO_U)$.

Let $\Art_k$ denote the category of local Artinian (commutative) $k$-algebras
with the residue field $k$. We are going to compare two deformation functors, \[ F_{1,n} ,
F_{\infty} : \Art_k \to \text{Sets}. \] Here $F_{1,n}(R)$ is the set of
isomorphism classes of families $C \to \Spec R$ with sections $p_1,\ldots, p_n$,
which reduce to $C_{1,n}$ upon the specialization $R\to k$. Note that this is nothing but the fibre of
$\wt{\UU}^{sns}_{1,n}(R) \to \wt{\UU}^{sns}_{1,n}(k)$ over the point corresponding to $C_{1,n}$.

Similarly, we define $F_\infty(R)$ as the fibre of $\MM_\infty(E_{1,n})(R)
\to \MM_\infty(E_{1,n})(k)$ over the class of the trivial $A_\infty$-structure. By Lemma 4.5.1 (i)
of \cite{P-ainf}, these correspond to equivalence classes of minimal $R$-linear
$A_\infty$-structures on $E_{1,n}$, such that upon specialization $R \to k$ we
get a formal $A_\infty$-algebra, i.e. the $A_\infty$-structure that is gauge equivalent to the trivial one.

Recall that the map $$\wt{\UU}^{sns}_{1,n}(k)\to
\MM_\infty(E_{1,n})(k)$$ (see \eqref{genus-1-a-inf-map})
sends the point corresponding to $C_{1,n}$ to the
class of the trivial $A_\infty$-structure (see Lemma \ref{Gm-action-compatibility-lem}(ii)).
Thus, we can consider the induced map 
\begin{equation}\label{def-functor-map}
F_{1,n}\to F_\infty 
\end{equation}
of deformation functors on $\Art_k$.  

\begin{prop} Assume that either $n\ge 3$, or $n=2$ and $\cha(k)\neq 2$ (resp., $n=1$ and $\cha(k)\neq 2,3$).
Then the morphism of deformation functors \eqref{def-functor-map} is an isomorphism, and the tangent space to $F_\infty$ is naturally isomorphic to \[ HH^2(E_{1,n})_{< 0} = HH^2(E_{1,n}) \]
\end{prop}
\Pf . The proof follows exactly the same line of argument as given in Prop.\ 4.5.4 of \cite{P-ainf},
using Corollary \ref{HH-C-1n-cor}, Lemma \ref{Ext-U-lem} and Lemma \ref{M-infty-repr-lem} (the latter is needed
to deduce that the functor $F_\infty$ is prorepresentable, hence, homogeneous). Note that, as we observed above,
the tangent space and obstruction space to $F_{1,n}$ can be identified with $\Ext^1_C(\bL_C,\OO(-D-p_i))$ and
$\Ext^2_C(\bL_C,\OO(-D-p_i))$, respectively.
\ed	

In particular, this leads to a computation of $HH^2(E_{1,n})$. Similarly to
\cite[Prop.\ 4.7.2]{P-ainf} we can identify $HH^2_{<0}(E_{1,n})$ with the
tangent space of the moduli scheme $\wt{\UU}^{sns}_{1,n}$ at zero. 
Furthermore, since \eqref{genus-1-a-inf-map} is compatible with the $\G_m$-action, this is a graded identification, hence, using Corollary \ref{tangent-cor} (resp. Theorem \ref{moduli-curves-thm} for $n\le 4$)
we get the ranks of $HH^2(E_{1,n})$. 

\begin{cor} \label{hh2} Over an arbitrary field $k$, we have for $n\geq 5$
\[ HH^2(E_{1,n}) = k^{(n-1)(n-2)/2}[1] \]
For $n=3,4$ we have:
\begin{align*}
	HH^2(E_{1,4}) &= k^3[1] \oplus k^2[2], \\ 
	HH^2(E_{1,3}) &= k[1] \oplus k^2[2] \oplus k[3]. 
\end{align*}
For a field $k$ with $\cha(k)\neq 2$,
\[ HH^2(E_{1,2}) = k[2] \oplus k[3] \oplus k[4]. \] 
For a field $k$ with $\cha(k) \neq 2 \text{\ or\ } 3$,
\[ HH^2(E_{1,1}) = k[4] \oplus k[6]. \]
\ed  \end{cor} 

\begin{rem}
It was shown in \cite{LP2} that
for a field $k$ with $\cha(k)=2$ one has
\[ HH^2(E_{1,1}) = k[1] \oplus k[3] \oplus k[4] \oplus k[6],   \]
while for a field $k$ with $\cha(k)=3$ one has
\[ HH^2(E_{1,1}) = k[2] \oplus k[4] \oplus k[6].    \]
Using the methods of this paper (or computing using an explicit resolution as in \cite{LP2})
one can show that for a field $k$ with $\cha(k) = 2$ one has
\[ HH^2(E_{1,2}) = k[1] \oplus k[2] \oplus k[3] \oplus k[4].\] 
\end{rem}

Finally, using the $\G_m$-action as in the proof of Theorem A of \cite{P-ainf}, we deduce our second main result.

\begin{thm}\label{ainf-thm} For $n\geq 3$, the map \eqref{genus-1-a-inf-map} induces an isomorphism of the moduli scheme $\wt{\UU}^{sns}_{1,n}$ with the moduli scheme of minimal $A_\infty$-structures
on $E_{1,n}$, up to a gauge equivalence. This isomorphism is compatible with the natural $\G_m$-actions. 
The same conclusion holds for $n=2$ and $\cha(k) \neq 2$ (resp., $n=1$ and $\cha(k)\neq 2,3$). 
\ed \end{thm} 

As a consequence, we get an interpretation of the moduli space $\ov{\MM}_{1,n}^\infty$ in terms of $A_\infty$-structures.

\begin{cor}\label{moduli} Under the assumptions of Theorem \ref{ainf-thm}, 
the quotient stack 
$$(\wt{\UU}^{sns}_{1,n} \backslash
	\{C_{1,n}\}) / \G_m =	\ov{\MM}_{1,n}^{\infty}$$ is isomorphic to the
	moduli stack of non-formal minimal $A_\infty$-structures on $E_{1,n}$, up to gauge equivalence and rescaling,
	or equivalently up to an $A_\infty$-equivalence of $A_\infty$-structures over $k^{n+1}$.
	\end{cor}

\Pf . We observe that all automorphisms of $E_{1,n}$ as an associative algebra over $k^{n+1}$ have form
$$A_i\mapsto \la_i A_i, \ B_i\mapsto \la\cdot\la_i^{-1}B_i,$$
for some invertible constants $\la,\la_i$. It is easy to check that the effect of such transformations on $\mu^n$
is exactly the rescaling $\mu^n\mapsto \la^{n-2}\mu^n$.
\ed

\begin{rem} 

Note that for every subset $S\sub\{1,\ldots,n\}$ we have a natural subquiver in
$Q_n$ such that the corresponding subalgebra is isomorphic to $E_{1,|S|}$. In
particular, we have $n$ subquivers $Q_{n-1}(i)\sub Q_n$ (where $i=1,\ldots,n$)
that give embeddings of $E_{1,{n-1}}$ into $E_{1,n}$.  Now given a minimal
$A_\infty$-structure $\mu^\bullet$ on $E_{1,n}$, for each $i$ we have a well
defined restriction $\mu^\bullet|_{Q_{n-1}(i)}$, which is a minimal
$A_\infty$-structure on $E_{1,{n-1}}$ (recall that we consider
$A_\infty$-structures that are unital with respect to the idempotents in $E_{1,n}$).
Therefore, we get maps \[ \tilde{\pi}_i : \MM_\infty(E_{1,n}) \to
\MM_\infty(E_{1,n-1}) \] for $i=1,\ldots,n$. Under the isomorphism
(\ref{genus-1-a-inf-map}), the map $\tilde{\pi}_i$ corresponds to forgetting the marked point
$p_i$ (over some open locus including the smooth curves). Thus, for $i=n$ this
morphism can be identified with the projection $U_n \to U_{n-1}$ of Proposition
\ref{better-equations-prop} (ii). \\

\end{rem}

Next, we observe that by Proposition \ref{perf}, we have that $HH^3(E_{1,n})
\simeq HH^3(C_{1.n})$ is finite-dimensional. Therefore, by \cite[Cor.
4.2.6]{P-ainf}, there is a natural isomorphism of functors
\[ \MM_\infty(E_{1,n}) \to \MM_d(E_{1,n}) \] for all $d
\geq N(n)$ for some sufficiently large $N(n)$, where the forgetful map is given by 
\[ (\mu^i)_{3 \leq i}  \to (\mu^i)_{3 \leq i \leq d} \]
In other words, the inverse limit
$\varprojlim_d \MM_d (E_{1,n})$ stabilizes. We next determine the
exact value of $N(n)$ for all $n$.

\begin{thm} Over an arbitrary field $k$, we have for $n \geq 4$, 
	\[ \MM_\infty(E_{1,n}) \simeq \MM_d (E_{1,n}) \text{\ \ for all $d \geq 4$. } \] 
For $n=3$, we have
	\[ \MM_\infty(E_{1,3}) \simeq \MM_d (E_{1,3}) \text{\ \ for all $d \geq 5$.  } \] 
For $n=2$, over a field $k$ with $\cha(k) \neq 2$, we have 
\[ \MM_\infty(E_{1,2}) \simeq \MM_d (E_{1,2}) \text{\ \ for all $d \geq 6$. } \]
For $n=1$, over a field $k$ with $\cha(k) \neq 2 \text{\ or\ } 3$, we have 
\[ \MM_\infty(E_{1,1}) \simeq \MM_d (E_{1,1}) \text{\ \ for all $d \geq 8$. } \] 
\end{thm}
\Pf .  Let $n\geq 5$. Abusing the notation we denote $U_n\times\Spec(k)$ simply as $U_n$.
By Corollary \ref{tangent-cor}, we know that
 \[ \Theta := \{a, (c_i, \ov{c}_i)_{4\le i\le n}, (c_{ij})_{4\le i<j\le n} \} \] is a set of
minimal generators of the algebra of functions on the affine scheme $U_n$,
where $a$, $c_i, \ov{c}_i$ and $c_{ij},\ov{c}_{ij}$ have degree $1$ with
respect to the natural $\G_m$ action. Furthermore, we have 
$k[U_n]=k[\Theta]/I$ where the ideal $I$ is generated by quadratic relations.
Now, the argument given in the proof of \cite[Prop. 4.7.2]{P-ainf} shows
that the algebra of functions on $\MM_{d}(E_{1,n})$ is isomorphic
\[ k[\Theta_{ \le {d-2}} ] / I_{\le {d-2}}, \]
where $\Theta_{\le d-2}$ (resp., $I_{\le d-2}$) is the set of elements in $\Theta$ (resp., $I$) of degree $\le d-2$. 
It follows
that $\MM_\infty(E_{1,n}) \simeq \MM_d (E_{1,n})$ for
$n\geq 5$ and $d \geq 4$. 

The proof in the remaining cases is similar. Again
we need to determine the algebra of functions on $U_n$ in these cases. These have been
worked out in Proposition \ref{better-equations-prop} in the cases $n=3,4$. 
The assertion follows as above, since $k[U_4]$ (resp., $k[U_3]$) is free with generators of degrees $\le 2$
(resp., $\le 3$).

The case of $n=2$ is worked out in Section \ref{case2} for $\cha(k) \neq 2$: we get 
that $k[U_2] = \A^3$ with generators $\alpha, \beta$ and
$\gamma$ of degrees $2,3$ and $4$, respectively. Hence, we have
$\MM_\infty(E_{1,2}) \simeq \MM_d (E_{1,2})$ for all
$d \geq 6$.

Finally, the case $n=1$ was worked out in \cite[Prop. 9]{LP2} for $\cha(k) \neq 2$ or $3$. (The general case is also studied in \cite{LP}). 
\ed

\begin{rem} We note that over a field $k$, the functor
	$\MM_3(E_{1,n})$ associating with $R$ the set of gauge
	equivalence classes of minimal $A_3$-structures on $E_{1,n}\otimes R$
	is represented by the affine space $HH^2(E_{1,n})_{-1}$ (see \cite[Thm.
	4.2.4]{P-ainf}). Assume that $n \geq 5$. By Corollary \ref{hh2}, it
	follows that \[ \MM_3(E_{1,n}) \simeq \A^{(n-1)(n-2)/2} \] 
		On the other hand, Theorem \ref{moduli-curves-thm}  together with Theorem
	\ref{ainf-thm} identifies the moduli scheme of $A_\infty$-structures on $E_{1,n}$ up to a gauge equivalence with the affine scheme 
	\begin{equation}\label{U-n-embedding-eq} 
	U_n \hra \A^{(n-1)(n-2)/2 }
	\end{equation} 
	In fact, this embedding can be identified with the natural map	\[ \MM_\infty (E_{1,n}) \hra \MM_3(E_{1,n}) \] 
	which sends $(\mu^i)_{i \geq 3}$ to $\mu^3$ by simply forgetting the higher products. 
	Finally, note that if we interpret $U_n\simeq\wt{\UU}^{sns}_{1,n}$ as moduli of curves
then we can view \eqref{U-n-embedding-eq} as an analog of the rational map from a $\G_m^g$-torsor over
$\MM_{g,g}$ to $\A^{g^2-g}$, defined in \cite{FisPol} in terms of triple products.
\end{rem} 

\end{section}

\end{document}